\newif\ifprint
\printfalse
\documentclass[twoside,10pt]{HandbookOfModuli}

\usepackage{amsmath,amsthm}
\usepackage{amssymb}
\usepackage{latexsym}
\usepackage[utf8]{inputenc}
\usepackage{textcomp}
\usepackage{color}
\usepackage[pdftex]{graphicx}
\usepackage{titlesec}
\usepackage{xspace}

\ifprint
	\input{giovanni} 
	\usepackage[final]{microtype}
\fi
\usepackage{bm}
\renewcommand{\mathbf}[1]{\bm{#1}} 

\ifprint
	\definecolor{linkred}{rgb}{0,0,0} 
	\definecolor{linkblue}{rgb}{0,0,0} 
\else
	\definecolor{linkred}{rgb}{0.7,0.2,0.2}
	\definecolor{linkblue}{rgb}{0,0.2,0.6}
\fi

\numberwithin{equation}{section} 

\usepackage[
	hyperindex,
	pagebackref,
	pdftex,
	pdftitle={Handbook of Moduli},
	pdfdisplaydoctitle,
	pdfpagemode=UseNone,
	breaklinks=true,
	extension=pdf,
	bookmarks=false,
	plainpages=false,
	colorlinks,
	linkcolor=linkblue,
	citecolor=linkred,
	urlcolor=linkred,
	pdfmenubar=true,
	pdftoolbar=true,
	pdfpagelabels,
	pdfpagelayout=TwoPage,
	pdfview=Fit,
	pdfstartview=Fit
]{hyperref}

\catcode`@=11 \baselineskip 4.5mm
\def\ps@handbook{\def\@oddhead{\hfill \leftmark \hfill\thepage }
\def\@evenhead{\thepage \hfill \rightmark \hfill}
\def\@oddfoot{}
\def\@evenfoot{}}
\def\@evenhead{}
\def\@oddfoot{}
\def\@evenfoot{\hfill\copyright\ China Higher Education Press}
\def\list#1#2{\ifnum \@listdepth >5\relax \@toodeep \else \global
\advance \@listdepth\@ne \fi \rightmargin \z@ \listparindent\z@
\itemindent\z@ \csname @list\romannumeral\the\@listdepth\endcsname
\def\@itemlabel{#1}\let\makelabel\@mklab \@nmbrlistfalse #2\relax
\@trivlist \parskip -\parsep \parindent\listparindent \advance
\linewidth -\rightmargin \advance\linewidth -\leftmargin \advance
\@totalleftmargin \leftmargin \parshape \@ne \@totalleftmargin
\linewidth \ignorespaces}
\catcode`@=12
\renewcommand\thesection{\arabic{section}}
\renewcommand\thesubsection{\arabic{subsection}}
\renewcommand\thesubsubsection{\arabic{subsubsection}}
\renewcommand{\theequation}{\thesection.\arabic{equation}}

\def\thebibliography#1{\section*{References}
\list{[\arabic{enumi}]}{\settowidth \labelwidth{[#1]} \leftmargin
\labelwidth \advance \leftmargin \labelsep \usecounter{enumi}}
\def\newblock{\hskip .11em plus .33em minus .07em} \sloppy
\clubpenalty 4000 \widowpenalty 4000 \sfcode`\.=1000 \relax}

\titleformat{\section}{\normalfont\large\bfseries}{\thesection.}{0.5em}{}
\titleformat{\subsection}{\normalfont\bfseries}{\thesection.\thesubsection.}{0.5em}{}
\titleformat{\subsubsection}[runin]{\normalfont\bfseries}{\thesection.\thesubsection.\thesubsubsection.}{0.5em}{}[\kern0.5em]

\usepackage{amssymb,latexsym,amsmath,amsthm,amsfonts,bm,hyperref}
\usepackage{pdfsync}
\input xy \xyoption {all}

\newtheorem{theorem}[equation]{Theorem}

\newtheorem{conjecture}[equation]{Conjecture}

\newtheorem{definition}[equation]{Definition}
\theoremstyle{remark}
\newtheorem{remark}[equation]{Remark}

\newenvironment{example}{\paragraph{\it Example.}}{\vskip0.4cm}

    \def\AC{{\mathcal{A}}}

    \def\MC{{\mathcal{M}}}
    \def\calN{{\mathcal{N}}}

\def\G{\Gamma}

\def\L{\Lambda}
\def\O{\Omega}

\newcommand{\nc}{\newcommand} \newcommand{\renc}{\renewcommand}

\newcommand{\rdots}{\mathinner{ \mkern1mu\raise1pt\hbox{.}
    \mkern2mu\raise4pt\hbox{.}
    \mkern2mu\raise7pt\vbox{\kern7pt\hbox{.}}\mkern1mu}}

\def\reg{{\mathrm{reg}}}

\def\red{{\mathrm{red}}}

\DeclareMathOperator{\Pic}{Pic}
\DeclareMathOperator{\Jac}{Jac}
\DeclareMathOperator{\Prym}{Prym}

\def\to{\rightarrow}

\def\longto{\longrightarrow}

\nc{\triright}{\stackrel{[1]}{\to}}
\nc{\longtriright}{\stackrel{[1]}{\longto}}

\nc{\Br}{\mathcal{B}}
\nc{\HotRR}{{}_R\mathcal{K}_R}
\nc{\HotR}{\mathcal{K}_R}
\nc{\excise}[1]{}
\nc{\defect}{\text{df}}
\nc{\h}[1]{\underline{H}_{#1}}

\nc{\Ga}{\mathbb{G}_a} 
\nc{\Gm}{\mathbb{G}_m} 

\nc{\Perv}{{\mathbf{P}}}

\nc{\IH}{{\mathrm{IH}}}

\nc{\ic}{\mathbf{IC}}

\nc{\gl}{{\mathfrak{gl}}}
\renc{\sl}{{\mathfrak{sl}}}
\renc{\sp}{{\mathfrak{sp}}}

\nc{\HBM}{H^{BM}}

\DeclareMathOperator{\End}{End} 

\DeclareMathOperator{\Spec}{Spec}

\DeclareMathOperator{\rk}{rk}

\newcommand{\bino}[2]{\mbox{ $#1 \choose #2$}}

\nc{\hp}{{\hat{p}}}
\renewcommand{\ell}{{\rm ell}}
\nc{\Gal}{\Gamma}
\nc{\st}{{\rm st}}
\nc{\cN}{{\check{\mathcal N}}}
\nc{\hN}{{\hat{\mathcal N}}}
\nc{\cM}{{\check{\mathcal M}}}
\nc{\hM}{{\hat{\mathcal M}}}
\nc{\N}{{\mathcal N}}
\nc{\cchi}{{\check{\chi}}}
\nc{\hchi}{{\hat{\chi}}}
\nc{\F}{{\mathbb F}}
\nc{\mF}{{\mathcal F}}
\nc{\mC}{{\mathcal C}}
\nc{\calM}{{\mathcal M}}
\nc{\calP}{{\mathcal P}}
\nc{\calT}{{\mathcal T}}
\nc{\Q}{{\mathbb Q}}
\nc{\J}{{\mathcal J}}
\nc{\higgsinfty}{${{\rm Higgs}_\infty}\ $}
\nc{\lrow}{\longrightarrow}
\nc{\Dol}{{\op{Dol}}}
\nc{\Hit}{{\op{Hit}}}
\nc{\DR}{{\op{DR}}}
\nc{\B}{{\op{B}}} 
\nc{\op}[1]{\mathop{\mathchoice{\mbox{\rm #1}}{\mbox{\rm #1}}
{\mbox{\rm \scriptsize #1}}{\mbox{\rm \tiny #1}}}\nolimits}
\nc{\down}[1]{{\phantom{riptstyle #1}
\hbox{$\left\downarrow\vbox to
    9.5pt{}\right.\nulldelimiterspace=0pt \mathsurround=0pt$}
\raisebox{.4ex}{$riptstyle #1$}}}
\nc{\Qu}{{\mathbb Q}}
\nc{\calC}{{\mathcal C}}
\nc{\calD}{{\mathcal D}}
\nc{\stack}[2]{
{\def\arraycolsep{0pt}\begin{array}{c}
riptstyle #1 \\ riptstyle #2 \end{array}} }

\nc{\Nm}{{\rm Nm}}
\nc{\rG}{{\rm G}}
\nc{\calE}{{\mathcal E}}
\nc{\calV}{{\mathcal V}}
\nc{\vv}{{\mathbf{v}}}
\nc{\w}{{\mathbf{w}}}
\nc{\calB}{{\mathcal B}}
\nc{\C}{{\mathbb C}}
\nc{\R}{{\mathbb R}}
\nc{\calK}{{\mathcal K}}
\nc{\calR}{{\mathcal R}}
\nc{\Est}{E_{\op{st}}}
\nc{\Hst}{H_{\op{st}}}
\nc{\Z}{{\mathbb Z}}
\nc{\Nat}{{\mathbb N}}
\nc{\E}{{\mathbb E}}
\nc{\U}{{\rm U}}
\nc{\dbar}{{\overline{\partial}}}
\nc{\K}{\mathbb K}
\nc{\Dirac}{{D\!\!\!\! \slash}} 
\nc{\Ha}{{\mathcal H}}
\nc{\Hy}{{\mathbb H}}
\nc{\GL}{{\rm GL}}
\nc{\SL}{{\rm SL}}
\nc{\PGL}{{\rm PGL}}
\nc{\V}{{\mathbb V}}
\nc{\M}{{\mathcal M}}
\nc{\I}{{\mathbb I}}
\renewcommand{\J}{{\mathbb J}}
\nc{\A}{{\mathcal A}}
\nc{\uPhi}{{\mathbf \Phi}}
\nc{\bM}{{\mathbf M}}
\nc{\calO}{{\mathcal O}}
\nc{\beq}[1]{\begin{eqnarray}\label{#1}}
\nc{\eeq}{\end{eqnarray}}
\nc{\bes}{\begin{eqnarray*}}
\nc{\ees}{\end{eqnarray*}}
\nc{\Left}[1]{\hbox{$\left#1\vbox to
    10.5pt{}\right.\nulldelimiterspace=0pt \mathsurround=0pt$}}
\nc{\Right}[1]{\hbox{$\left.\vbox to
    10.5pt{}\right#1\nulldelimiterspace=0pt \mathsurround=0pt$}}

\begin{document}
\setcounter{page}{1}
%
%
\long\def\replace#1{#1}

%
%
{\title{Global topology of the Hitchin system} }
%
%
\author{Tam\'as Hausel}
\address{24-29 St Giles', Mathematical Institute,  Oxford, OX1 3LB, United Kingdom}
\email{{hausel@maths.ox.ac.uk}}

%
%

\begin{abstract}
	
{Here we survey several results and conjectures on the cohomology of the total space of the Hitchin system:
	the moduli space of semi-stable rank $n$ and degree $d$ Higgs bundles on a complex algebraic curve $C$. 
	The picture emerging is a dynamic mixture of ideas originating in theoretical physics such as gauge theory and mirror symmetry, Weil conjectures in arithmetic algebraic geometry,
	representation theory of finite groups of Lie type and Langlands duality in number theory. }

\end{abstract}  

\maketitle
\thispagestyle{empty}
%
%

\section{Introduction}
\label{intro}

Studying the topology of moduli spaces in algebraic geometry could be considered the first approximation of understanding the moduli problem. We start with an example which is one of the original examples of moduli spaces constructed by Mumford \cite{mumford} in 1962 using Geometric Invariant Theory. Let $\calN$ denote the moduli space of semi-stable rank $n$
degree $d$ stable bundles on a smooth complex algebraic curve $C$; which turns out to be a projective variety.  In 1965 Narasimhan and Seshadri \cite{narasimhan-seshadri} proved that the space $\N$ is canonically diffeomorphic
with the manifold $$
{\calN}_\B:= \{ A_1,B_1,\dots, A_g,B_g \in {\rm U}(n)\ | 
		 [A_1, B_1]\dots [A_g, B_g] = \zeta_n^d I \}/{\rm U}(n),
$$
of twisted representations of the fundamental group of $C$ to $\U(n)$, where $\zeta_n=\exp({2\pi i/n})$. Newstead \cite{newstead} in 1966 used this latter description to determine the Betti numbers of
$\calN_\B$ when $n=2$ and $d=1$. Using ideas from algebraic number theory for function fields Harder \cite{harder}  in 1970 counted the rational points $\calN(\F_q)$ over a finite field when $n=2$  and $d=1$ and  compared his formulae to Newstead's results to find
that the analogue of the Riemann hypothesis, the last remaining  Weil conjectures, holds in this case. By the time Harder and Narasimhan  \cite{harder-narasimhan} in 1975 managed to count $\#\calN(\F_q)$ for any $n$ and coprime $d$ the last of the Weil conjectures had been proved
by Deligne \cite{deligne} in 1974, and their result in turn yielded \cite{desale-ramanan} recursive formulae for the Betti numbers of $\calN(\C)$ .

The same recursive formulae were found using a completely different method by Atiyah and Bott \cite{atiyah-bott} in 1981. They studied the topology of $\calN$ in another reincarnation with origins in theoretical physics. The Yang-Mills equations on $\R^4$
appeared in theoretical physics as certain non-abelian generalization to Maxwell's equation; and were used to describe aspects of the standard model of the physics of elementary particles. Atiyah and Bott considered the analogue of these equations in 2 dimensions, more precisely
 the solution  space $\N_{\op{YM}}$ 
 of Yang-Mills connections on a  differentiable Hermitian vector bundle of rank $n$ and degree $d$ on the Riemann surface $C$ modulo gauge transformations. Once again we have canonical diffeomorphisms  $\N_{\op{YM}}\cong \calN_\B \cong \calN$. 
Atiyah--Bott \cite{atiyah-bott} used this gauge theoretical approach to study
the topology of  
$\calN_{\op{YM}}$ obtaining in particular recursive formulae for the Betti numbers which turned out to be essentially the same as those arising from Harder-Narasimhan's arithmetic approach. Atiyah-Bott's approach has been greatly generalized
in Kirwan's work \cite{kirwan} to study the cohomology of K\"ahler quotients in differential geometry and the closely related GIT quotients in algebraic geometry. By now we have a fairly comprehensive understanding of the cohomology (besides
Betti numbers: the ring structure, torsion, K-theory etc) of $\N$ or more generally of compact K\"ahler and GIT quotients.

In this survey we will be interested in a hyperk\"ahler analogue of the above questions. It started with Hitchin's seminal paper \cite{hitchin} in 1987. Hitchin \cite{hitchin} investigated the two-dimensional reduction of the four dimensional Yang-Mills equations
over a Riemann surface, what we now call Hitchin's self-duality equations on a rank $n$ and degree $d$ bundle on the Riemann surface $C$. Using gauge theoretical methods one can construct  $\M_\Hit$ the space of solutions to Hitchin's self-duality 
equations modulo gauge transformations. $\M_\Hit$ can be considered as the hyperk\"ahler analogue of $\N_{\op{YM}}$. When $n$ and $d$ are coprime $\M_\Hit$ is a smooth, albeit non-compact, differentiable manifold with a natural hyperk\"ahler metric. The latter means a metric with three K\"ahler forms corresponding to complex structures $I,J$ and $K$, satisfiying quaternionic relations. Hitchin \cite{hitchin} found that the corresponding algebraic geometric moduli space is the moduli space $\M_{\Dol}$ of semi-stable Higgs bundles $(E,\phi)$ of rank $n$ and degree $d$  and again we have the natural diffeomorphism $\M_{\Hit}\cong \M_{\Dol}$ which in fact is complex analytical in the complex structure $I$ of the hyperk\"ahler metric on $\M_{\Hit}$. $\M_\Dol$ is the hyperk\"ahler analogue of $\N$. Using a natural Morse function on $\M_\Hit$ Hitchin was able to determine the Betti numbers of $\M_\Hit$ when $n=2$ and $d=1$, and this work was extended by Gothen \cite{gothen} in 1994 to the $n=3$ case. Hitherto the Betti numbers for $n>3$
have  not been found, but see \S\ref{conclusion}. 

In this paper we will survey several approaches to get cohomological information on $\M_\Hit$. The main difficulty lies in the fact that $\M_\Hit$ is non-compact. However as Hitchin \cite{hitchin} says 
\vskip.2cm
"... the moduli space of all solutions turns out to be a
manifold with an extremely rich geometric structure". 
\vskip.2cm
\noindent Due to this surprisingly rich geometrical structures on $\M_\Hit$
we will have several different approaches to study $H^*(\M_\Hit)$, some motivated by ideas in theoretical physics, some by arithmetic and some by Langlands duality. 

As mentioned $\M_\Hit$ has a natural hyperk\"ahler metric and in complex structure $J$ it turns out to be complex analytically isomorphic to the character variety
$$
\M_\B:= \left\{ \right. A_1,B_1,\dots, A_g,B_g \in \GL_n(\C)\ | 
		[A_1, B_1]\dots [A_g, B_g] =\zeta_n^d I \left. \right\}/\!/\GL_n(\C)
$$
of twisted representations of the fundamental group of $C$ to $\GL_n$. The new phenomenon here is that this is a complex variety. Thus cohomological 
information could be gained by counting points of it over finite fields.  Using the character table of the finite group $\GL_n(\F_q)$ of Lie type this was
accomplished in \cite{hausel-villegas}. The calculation leads to some interesting conjectures on the Betti numbers of $\M_\B$. This approach was surveyed
in \cite{hausel-mln}. We will also mention this approach in \S\ref{betti} but here will focus on its connections to one more aspect of the geometry of $\M_\Hit$. Namely, that it carries the structure of an integrable system.

This integrable system, called the Hitchin system, was defined in \cite{Hitchin} by taking the characteristic polynomial of the Higgs field, and has several remarkable properties. First of all
it is a completely integrable system with respect to the holomorphic symplectic structure on $\M_\Dol$ arising from the hyperk\"ahler
metric. Second it is proper, in particular the generic fibers are Abelian varieties. 
 In this paper we will concentrate on the topological aspects of the Hitchin system on $\M_\Dol$. 

The recent renewed interest in the Hitchin system could be traced back to two major advances. One is Kapustin-Witten's \cite{kapustin-witten}  proposal in 2006, that S-duality  gives a physical framework for the geometric Langlands correspondence.
In fact, \cite{kapustin-witten} argues that this S-duality reduces to mirror symmetry for Hitchin systems for Langlands dual groups. This point of view has been expanded and explained
in several papers such as \cite{gukov-witten}, \cite{frenkel-witten}, \cite{witten}.  One of our main
motivation to study $H^*(\M_\Dol)$ is to try to prove the topological mirror symmetry proposal of \cite{hausel-thaddeus},
which is a certain agreement of Hodge numbers of the mirror Hitchin systems; which we now understand as a cohomological shadow of Kapustin-Witten's S-duality (see \S\ref{solving}.\ref{shadow}). This topological mirror symmetry proposal is mathematically well-defined
and could be tested by using the existing techniques of studying the cohomology of $H^*(\M_\Dol)$. Indeed it is a theorem
when $n=2,3$ using \cite{hitchin} and \cite{gothen}. Below we will introduce the necessary formalism to discuss these
conjectures and results.

Most recently Ng\^{o} proved \cite{ngo} the fundamental lemma in the Langlands program in 2008  by a detailed study of the topology
of the Hitchin map over a large open subset of its image.   Surprisingly, Ng\^{o}'s geometric approach will be intimately
related to our considerations, even though Ng\^{o}'s main application is the study of the cohomology of a  {\em singular fiber} of the Hitchin map.
We will explain at the end of this survey some connections of Ng\^{o}'s work to our studies of the global topology of the Hitchin map. 

We should also mention that both Kapustin-Witten and Ng\^{o} study Hitchin systems for general reductive groups, for simplicity we will concentrate on the groups $\GL_n,\SL_n$ and $\PGL_n$ in this paper.

Our studies will lead us into a circle of ideas relating arithmetic and the Langlands program to the physical ideas from gauge theory,  S-duality and mirror symmetry in the study of the global topology of the Hitchin system. This could be considered the hyperk\"ahler analogue 
of the fascinating parallels\footnote{For  a recent survey on these parallels see \cite{kirwan-etal}.} between the arithmetic approach of Harder--Narasimhan and the gauge theoretical approach of Atiyah--Bott in the study of $H^*(\calN)$.
\vskip.2cm
\begin{paragraph}{\bf Acknowledgement} The basis of this survey are the notes taken by Gergely B\'erczi, Michael Gr\"ochenig
	and Geordie Williamson from a preparatory lecture course "Mirror symmetry, Langlands duality and the Hitchin system" the author gave in Oxford in Hilary 
	term 2010. I would like to thank them and the organizers of several lecture courses where various subsets of these notes were delivered in 2010. These include the Spring school on  "Geometric Langlands and gauge theory"
	at CRM Barcelona, a lecture series  at the Geometry and Quantum Field Theory cluster at the University of Amsterdam, the "Summer School on the Hitchin fibration" at the University of Bonn and  the Simons Lecture Series at Stony Brook. I would also like to thank Mark de Cataldo, Iain Gordon, Jochen Heinloth, Daniel Huybrechts, Luca Migliorini, David Nadler, Tony Pantev and especially an anonymous referee for helpful comments. The author was supported by a Royal Society
	University Research Fellowship. 
\end{paragraph}

\subsection{Mirror symmetry}\label{mirror} In order to motivate the topological mirror symmetry conjectures of \cite{hausel-thaddeus}, we give some background information on three of the main topics involved starting with mirror symmetry.

 Considerations of mirror symmetry appeared  in various forms in string theory at the end of the 1980's.   It entered into the realm of mathematics via the work \cite{candelas-etal} of Candelas-de\! la\! Ossa-Green-Parkes in 1991   by formulating mathematically precise (and surprising) conjectures on the number of rational curves
in certain Calabi-Yau 3-folds. 

		 Mathematically,  mirror symmetry relates the symplectic geometry of a Calabi-Yau manifold
		$X$ to the complex geometry of its mirror Calabi-Yau $Y$ of the same dimension. Such an unexpected duality  between two previously separately studied 
		fields in geometry  has caught the
		interest of several mathematicians working in related fields.  The literature of mirror symmetry is vast, here we only mention monographs and conference proceedings on the subject such as 
		\cite{mirror1,mirror2,mirror3}.
		
		 	First aspect of checking a mirror symmetry proposal is the {\em topological mirror test} \beq{otms}h^{p,q}(X)=h^{\dim Y-p,q}(Y)\eeq that the Hodge diamonds have to be mirror to each other. If we introduce the notation
		  $$E(X;x,y)=\sum (-1)^{i+j} h^{i,j}(X)x^iy^j$$  the symmetry translates to \beq{functional} E(X;x,y)=x^{\dim Y
		 } E\left(Y;1/x,y\right).\eeq As we will be interested in the case when  the Calabi-Yau manifolds are
		hyperk\"ahler
we note that  a {\em compact} hyperk\"ahler manifold $X$ satisfies $h^{p,q}(X)=h^{\dim X-p,q}(X)$ and so the above topological mirror test 
	simplifies to the agreement \beq{ctms}  h^{p,q}(X)=h^{p,q}(Y)\eeq of Hodge numbers when both $X$ and $Y$ are hyperk\"ahler. 	
	  Kontsevich \cite{kontsevich} in 1994 suggested that mirror symmetry is underlined by a more fundamental  {\em homological mirror symmetry},  which identifies two derived categories   $${\mathcal D}^b(Fuk(X,\omega))\sim {\mathcal D}^b(Coh(Y,I))$$  the Fukaya category 		${\mathcal D}^b(Fuk(X,\omega))$ of certain decorated Lagrangian subvarieties of the symplectic manifold $(X,\omega)$ and the derived category
		of coherent sheaves ${\mathcal D}^b(Coh(Y,I))$ on the mirror $Y$, considered as a complex variety. This suggestion has been checked in several examples 
		\cite{polischuk-zaslow,seidel} and has been the starting point of a large body of mathematical research.
		 
A more geometrical proposal is contained in the work of Strominger-Yau-Zaslow \cite{strominger-etal} in 1996. They  suggested a geometrical construction how to obtain the mirror $Y$ from any given Calabi-Yau $X$. It suggested that $Y$ can be constructed as the moduli space of certain special Lagrangian submanifolds of  $X$ together with a flat line bundle on it; the picture arising then can be described as

		
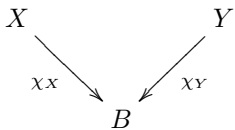
\[
		\xymatrix{
		X \ar[dr]_{{\chi}_X} & & Y \ar[dl]^{{\chi}_Y} \\
		& B
		}
		\]
	where $\chi_X$ and $\chi_Y$ are special Lagrangian fibrations on $X$ and $Y$ respectively, with generic fibers
	dual middle-dimensional tori.
	Until this  \cite{strominger-etal}  proposal there was no general geometrical conjectures how one might construct $Y$ from $X$. Even though there have been a fair amount of work today we do not have an example where the geometrical proposal of \cite{strominger-etal} is completely implemented for the original case of Calabi-Yau $3$-folds. But see \cite{gross} for a survey of recent developments. 
		 
		We conclude by noting that many mathematical predictions of mirror symmetry have been confirmed but we still have no general understanding yet.

\subsection{Langlands duality}
 Here we give some technically simplified remarks as a sketch of a basic introduction to a few ideas in the Langlands program in number theory; more details and references can be found in the  surveys \cite{frenkel,gelbart}. 
 
 As a first approximation,  the Langlands correspondence aims to describe a central object in number theory: the absolute Galois group ${\rm Gal}(\overline{\Q}/\Q)$ via representation theory. 
More precisely	let $\rG$  be a reductive group (over the complex numbers these are just the complexifications of compact Lie groups) and $^L\rG$ its Langlands dual, which one can obtain by dualizing the root datum of $\rG$.
	For example for the groups of our concern in the present paper:  $^L\GL_n=\GL_n$; $^L\SL_n=\PGL_n$, $^L\PGL_n=\SL_n$.
	 Langlands in 1967 conjectured that one can find a correspondence between the set isomorphism classes 
	of certain continuous homomorphisms 
${\rm Gal}(\overline{\Q}/\Q)\! \to\! \rG(\overline{\Q}_l)$ (for $\rG=\GL_n$ these are just  $n$-dimensional representations) and  	 automorphic (certain infinite dimensional) representations  of $^L\rG({\mathbb A}_\Q) $ over the ring of ad\`eles ${\mathbb A}_\Q$. 
The motivation for understanding the representations of ${\rm Gal}(\overline{\Q}/\Q)$ is that it describes the absolute Galois group itself via the Tannakian formalism. 

	Langlands built his programme as a non-abelian generalization of the ad\`elic description of class field theory 
	which can be understood as the 
	$\rG=\GL_1$  case of the above. Indeed representations of  ${\rm Gal}(\overline{\Q}/\Q)$ to $\GL_1$ describe the abelianization ${\rm Gal}_{ab}(\overline{\Q}/\Q)$ which describes all finite abelian extensions of $\Q$. 
	In the case of $\rG=\GL_2$ elliptic curves enter naturally, via the action of the absolute Galois group on their two-dimensional first \'etale cohomology. The corresponding objects on the automorphic side are modular forms, and the Langlands program in this case can be seen to reduce to the Shimura-Taniyama-Weil conjecture which is now a theorem due to the work of Wiles and others. 
	
 	The ad\`elic description of class field theory was partly motivated by the deep analogy between $\Q$,  or more generally number fields (i.e. finite extensions of the rationals), and the function field $\F_q(C)$, where $C/\F_q$ is an algebraic curve. This analogy proved
powerful in attacking problems in algebraic number theory and in the Langlands program specifically. Many of the conjectured properties of number fields can be formulated for function fields, where one can use the techniques of algebraic geometry and succeed even in situations where the number field case remains open. Such an example is the Riemann hypothesis, which is still open over $\Q$, but was proved  over function fields $\F_q(C)$ of curves by Weil \cite{weil} in 1941 and for any algebraic variety by Deligne \cite{deligne} in 1974. 
 
 For an example of the more harmonious interaction between the number field and function field side of algebraic number theory we mention that  
Ng\^{o}'s recent proof \cite{ngo} of the  fundamental lemma for $\F_q(C)$ (where he used the topology of the Hitchin system, for more details see \S\ref{solving}.\ref{TMStoFL})
yielded the fundamental lemma for the number field case, due to previous work by Waldspurger \cite{waldspurger}. 
 
If we replace $\F_q$ by $\C$, i.e. consider the complex function field $\C(C)$ for complex curves $C/\C$, i.e. Riemann surfaces,  then we lose some of the analogies, but still
one can reformulate a version of the Langlands correspondence mentioned above. This is the {\em Geometric Langlands Correspondence} as was proposed by \cite{laumon,beilinson-drinfeld}. Originally it  was 
a correspondence between isomorphism classes of $\rG$-local systems on $C$ (analogues of the representations of the absolute Galois group) and certain Hecke eigensheaves (the analogues of automorphic representations) on the stack ${\rm Bun}_{^L\rG}$ of $^L\rG$-bundles on $X$. In these lectures our point of view will be that as a cohomological shadow of these considerations 
of the Geometric Langlands programme one can extract agreements of certain Hodge numbers 
of Hitchin systems for Langlands dual groups.

\subsection{Hitchin system}\label{integrablesystem}   Here we collect some of the basic ideas of completely integrable Hamiltonian systems and comment on the history of a large class of them: the Hitchin systems. An extensive account for the former can be found in \cite{arnold} while \cite{donagi-markman} details the latter.

Recall that a Hamiltonian system is given by a symplectic manifold $(X^{2d},\omega)$ and a Hamiltonian - or energy- function $H:X\to \R$. The corresponding 
	$X_H$ Hamiltonian vector field  is defined by the property $$dH=\omega(X_H,.).$$ The dynamics of the Hamiltonian system is given by the flow - the one parameter group
	of diffeomorphisms - generated by $X_H$. A function
	  $f:X\to \R$ is called a {\em first integral} if $$X_H f=\omega(X_f,X_H)=0$$ holds. The condition is equivalent to say that $f$ is constant along the flow of the system. 
	Note that  $\omega(X_H,X_H)=0$ as $\omega$ is alternating and so $H$ is constant along the flow - which is sometimes referred to as the law of conservation of energy.

	  We say that the Hamiltonian system is {\em completely integrable} if there is a map $$f=(H=f_1,\dots,f_d):X\to \R^d,$$ which is generic (meaning that  $f$ is generically a submersion) such that $\omega(X_{f_i},X_{f_j})=0$. 
	  The generic fibre of $f$ then has an action of $\R^d=\langle X_{f_1},\dots,X_{f_d}\rangle$  and so when $f$ is proper the generic fibre can be identified with a torus $(S^1)^d$. In such
	cases one has a fairly good control of the dynamics of the system (on the generic fiber it is just an affine motion) hence the name integrable. 
Several  examples arise in classical mechanics such as the Euler and Kovalevskaya tops and the spherical pendulum. For more details see \cite{arnold}.

Here	 we will be concerned  with the complexified or algebraic version of integrable Hamiltonian systems. Thus we consider a complex $2d$ dimensional manifold $X$ with a holomorphic symplectic
	$2$-form. $f:X\to \C^d$ now is a  algebraically completely integrable system when $f$ is generically a submersion and $\omega(X_{f_i},X_{f_j})=0$. If $f$ turns out to be proper, the generic fiber then will become a torus with a complex structure - in the algebraic case-  an Abelian variety.  This is the case for a large class of examples:  the {\em Hitchin systems}. As we will see in more detail below the Hitchin system is attached to the cotangent bundle of the moduli space  of stable $\rG$-bundles on a complex curve $X$. 
Originally  it appeared in Hitchin's study \cite{hitchin} 
 of the 2-dimensional reduction of the Yang-Mills equations from $4$-dimensions. Here we will follow a more algebraic approach.

\section{Higgs bundles and the Hitchin system}
\label{basic}
\subsection{The moduli space of vector bundles on a curve} 
Let $C$ be a complex projective
curve of genus $g > 1$. We fix integers $n > 0$ and $d \in
\mathbb{Z}$. We asssume throughout that $(d,n) = 1$.
Using Geometric Invariant Theory \cite{mumford,newstead2} one can construct
$\mathcal{N}^d $; the moduli space of isomorphism classes of 
stable rank $n$ degree $d$ vector bundles (equivalently $\GL_n$-bundles)
    on $C$.

We recall that a vector bundle is called \emph{semistable} if every subbundle $F$ satisfies
\[
\mu(F) = \frac{\deg F}{\rk F} \leq \mu(E) = \frac{\deg E }{\rk E}
\]
A vector bundle is \emph{stable} if one has strict inequality above
for all proper subbundles.

In general one has to be careful in constructing such
moduli spaces as special care has to be taken for the non-trivial automorphisms of strictly semi-stable objects. However, as we assume $(d,n) = 1$ the notions of
semi-stability and stability clearly agree. In particular, we can 
 conclude that $\mathcal{N}^d$ is   smooth and projective of dimension $d_n=n^2(g-1)+1$.

Consider the determinant morphism
\[
\det: \mathcal{N}^d \to \Jac^d(C)
\]
which sends a vector bundle of rank $n$ to its highest exterior power
$\Lambda^{n} E$. Choose $\Lambda \in \Jac^d(C)$ and define
$\check{\mathcal{N}}^{\Lambda} := \det^{-1}(\Lambda)$. When $\Lambda=\calO_C$ is the trivial bundle, 
the vector bundles
in $\check{\mathcal{N}}^{\Lambda}$ are exactly the $\SL_n$-bundles, for general $\Lambda$ we can think of
points in $\check{\mathcal{N}}^{\Lambda}$ as ``twisted
$\SL_n$-bundles''. Tensoring with an $n$th root of $\Lambda_1\Lambda_2^{-1}$ gives
an isomoprhism $\check{\calN}^{\Lambda_1}\cong\check{\calN}^{\Lambda_2}$ thus the
isomorphism class of
$\check{\mathcal{N}}^{\Lambda}$ does not depend  on
the choice of $\Lambda \in \Jac^d(C)$. We often abuse notation and
write $\check{\mathcal{N}}^d$ instead of
$\check{\mathcal{N}}^{\Lambda}$.

The abelian variety $\Pic^0(C) = \Jac^0(C)$ acts on $\mathcal{N}^d$
via
\[
(L,E) \mapsto L \otimes E.
\] As Seshadri showed in \cite{seshadri} the quotient of a normal variety by an Abelian
variety always exist, so we can define the moduli space of degree $d$ $\PGL_n$ bundles:
\[
\hat{\mathcal{N}}^d := \mathcal{N}^d/\Pic^0(C).
\] 
The embedding $\cN^\Lambda\subset \calN^d$ induces $$\hat{\mathcal{N}}^d \cong
\check{\mathcal{N}}^{\Lambda}/\Gamma.$$ Here $\Gamma := \Pic^0(C)[n]$ is
 the group of $n$-torsion points of the Jacobian, isomorphic to the finite group $\Z_n^{2g}$.  It acts on $\check{\mathcal{N}}^{\Lambda}$ by tensorization. Hence
$\hat{\mathcal{N}}^d$ is a projective orbifold.

\subsubsection{Cohomology of moduli spaces of bundles}

The cohomology\footnote{In this paper cohomology is with rational coefficients; unless indicated otherwise.} of $\mathcal{N}^d$, $\check{\mathcal{N}}^d$ and
$\hat{\mathcal{N}}^d$ is well understood. The structure of the cohomology
rings can be described by finding universal generators and all the relations
between them \cite{atiyah-bott, thaddeus,   earl-kirwan}.

Here we only comment on the additive structure in more detail. 
In 1975  Harder and Narasimhan \cite{harder-narasimhan} 
 obtained recursive formulae for the number of points of these
varieties over finite fields. It is then possible to use the Weil
conjectures (which had been
proven the year before by Deligne \cite{deligne}) to obtain formulae for the
Betti numbers. In 1981 Atiyah and Bott gave a different 
gauge-theoretic proof \cite{atiyah-bott}.

The main application in Harder and Narasimhan's paper is the following:

\begin{theorem}[\cite{harder-narasimhan}] \label{harnar}
The finite group $\Gamma$ acts trivially on
$H^*(\check{\mathcal{N}}^d)$. In particular, we have $H^*(\check{\mathcal{N}}^d) = H^*(\hat{\mathcal{N}}^d)$.
\end{theorem}
\begin{remark}
This result is difficult to prove and relies on showing that the
varieties $\check{\mathcal{N}}^d$ and $\hat{\mathcal{N}}^d$ have the
same number of points over finite fields. The analogue of this result
is false in the context of the moduli space of Higgs bundles as was
already observed by Hitchin in \cite{hitchin} for $n=2$. Interestingly for us, this will lead to the non-triviality of our topological mirror tests.

\end{remark}

\subsection{The Hitchin system}
We now consider the Hitchin system, which will be an integrable system on the
cotangent bundle to the moduli spaces considered in the previous
section. As in the previous section, fix $n$ and $d$ and abbreviate
$\calN := \calN^d$.

The cotangent bundle $T^*\mathcal{N}$ is an algebraic
variety. The ring of regular functions
$\mathbb{C}[T^*\mathcal{N}]$ turns out to be finitely-generated as will be proven below. 
The \emph{Hitchin map} then is simply the affinization:
\beq{hitchinfirst}
\chi: T^*\mathcal{N} \to \mathcal{A} = \Spec(\mathbb{C}[T^*\mathcal{N}]).
\eeq

We now describe this map more explicitly. For a point $[E] \in
\mathcal{N}$ standard deformation theory 
gives us an identification
\[ T_{[E]}\mathcal{N} = H^1(C, \End(E)).
\]
Applying Serre duality we obtain
\[
T_{[E]}^*\mathcal{N} = H^0(C, \End(E) \otimes K),\]
where $K$ denotes the canonical bundle of $C$. An element
\[ \phi \in H^0(C, \End(E) \otimes K)\]
is called a \emph{Higgs field}. Morally, it
can be thought of as a matrix of one-forms on the curve.
As such if we consider the characteristic polynomial of $(E,\phi) \in T^*\mathcal{N}$ then it will have the form
\[
t^n + a_1 t^{n-1} + \dots + a_n
\]
where $a_i \in H^0(K^i)$. For example $a_n \in H^0(K^n)$ is the
determinant of the Higgs field.

As we will prove below the Hitchin map \eqref{hitchinfirst} then has the explicit description
\begin{eqnarray*}
\chi : T^*\mathcal{N} \to & \AC := \bigoplus_{i = 1}^n H^0(K^i) \\
(E, \phi) \mapsto & (a_1, a_2, \dots, a_n)
\end{eqnarray*}
The affine space $\AC$ is called the \emph{Hitchin base}. 

In the $\SL_n$-case we have
\[
T_{[E]}^*\check{\mathcal{N}}^{d} = H^0(\End_0(E) \otimes K)\]
that is, a covector at $E$ is given by a \emph{trace free}
Higgs field. Thus in this case the Hitchin base is 
\[
 \A^0:=\bigoplus_{i = 2}^{n}H^0(C,K^i).
\] and the Hitchin map \beq{cchi} \cchi:T^*\check{\mathcal{N}}^{d}\to \A^0.\eeq

As the characteristic polynomial of the Higgs field does not change when the Higgs bundle is tensored with a line bundle, the action
of $\G$ on $T^*\cN$ is along the fibers of $\cchi$ and so $\cchi$ descends to the quotient
which gives the $\PGL_n$ Hitchin map:
$$\hchi:(T^*\cN)/\Gamma\to \hat{\A}=\A^0.$$

Recall that $T^* \calN$ is an algebraic symplectic variety with the canonical Liouville symplectic form. 

\begin{theorem}[Hitchin, 1987]
If $\chi_i$, $\chi_j\in \C[T^*{\mathcal N}]$ are two coordinate functions, then they Poisson commute, i.e. $\omega(X_{\chi_i},X_{
\chi_j}) = 0$. We have $\dim( \mathcal{A}) = \dim( \mathcal{N})$ and  the generic fibres of $\chi$ are open subsets of abelian varieties. Therefore we have an algebraically completely integrable Hamiltonian system.
\end{theorem}

As a next step we will projectivize  the Hitchin map $\chi: T^*\mathcal{N} \to
\mathcal{A}$.
 Recall that a complex point in $T^* \mathcal{N}$ is
given by a pair $(E,\phi)$. In order to projectivize we need to allow $E$ to become unstable.

\begin{definition}
A \emph{Higgs bundle} is a pair $(E,\phi)$ where $E$ is a vector bundle and
$\phi \in H^0(C, \End(E) \otimes K)$ is a Higgs field.
\end{definition}
The definition for semi-stability and stability for Higgs-bundles is
almost the same as for vector bundles except we only consider
$\phi$-invariant subbundles. The moduli-space of semi-stable Higgs
bundles is denoted by $\mathcal{M}_\Dol^d$ and often abbreviated as $\M^d$ 
or even $\M$. It  is a non-singular
quasi-projective variety, having $T^*\mathcal{N}$ as an open
subvariety.

It is straightforward to extend $\chi: \mathcal{M}^d \to \mathcal{A}$. The
following result shows that we have succeeded in projectivizing the
Hitchin map:

\begin{theorem}[Hitchin 1987 \cite{hitchin}, Nitsure 1991 \cite{nitsure}, Faltings 1993 \cite{Faltings}] \label{proper}
$\chi$ is a proper algebraically completely integrable Hamiltonian system. Its generic fibres are abelian varieties.
\end{theorem}

\begin{remark} Note that $T^*{\mathcal{N}}\subset {\mathcal{M}}$ is open dense with
	complement of codimension greater than $1$ \cite[Proposition 6.1.iv]{Hitchin}. Thus $\C[T^*{\mathcal{N}}]\cong \C[\mathcal{M}]$. Therefore
	 the affinization of $\mathcal{M}$ (which must be the Hitchin map as it
	is proper to an affine space) restricts to the affinization of $T^*{\mathcal{N}}$, which
	justifies our unorthodox introduction of the Hitchin map above. 
	\end{remark}
	
	\subsection{Hitchin systems for $\SL_n$ and $\PGL_n$}

It is straightforward now to compactify the $\SL_n$-Hitchin map $T^*\cN^d\to \A^0$. We consider $\cM_\Dol^\Lambda$ the moduli space of isomorphism classes of semi-stable  Higgs bundle $(E,\phi)$ of rank $n$, $\det E = \Lambda$ and trace-free $$\phi \in H^0(\End_0(E) \otimes K)$$ Higgs field.  Again the isomorphism class of $\cM_\Dol^\Lambda$ only depends on $d$, so we will simplify our notation to $\cM_\Dol^d$, $\cM^d$ or even $\cM$ for the $\SL_n$ Higgs moduli space.

As in the $\GL_n$-case the Hitchin map $\check{\chi}:\check{\M}^d\to \A^0$ is given by the coefficients of the characteristic polynomial, it is proper and a completely integrable system. We also have that $T^*\check{\mathcal{N}}^d \subset   \check{\mathcal{M}}^d$ open and dense.

Let us recall the two constructions of the moduli space of $\PGL_n$-Higgs bundles. The cotangent bundle $T^* \Pic^0(C) = \Pic^{0}(C) \times H^0(C,K)$ is a group. It acts on $\mathcal{M}^d$ by
\[
	(L,\varphi) \cdot (E,\phi) \mapsto (L\otimes E, \varphi + \phi)
\]
This induces an action of $\Gamma = \Pic^0[n]$ on $\check{\mathcal{M}}^d$. Then we may either define the $\PGL_n$-moduli space as
\[
	\hat{\mathcal{M}}^d = \mathcal{M}^d / T^* \Pic^0(C) \cong \chi^{-1}(\mathcal{A}^0) / \Pic^0(C)
\]
or equivalently as the orbifold
\[
	\hat{\mathcal{M}}^d=\check{\mathcal{M}}^d/\Gamma.
\]
The second quotient tells us that we obtain an orbifold. Since $\check{\chi}$ is compatible with the $\Gamma$ action, we obtain a well-defined proper Hitchin map
\[
	\hat{\chi}:\; \hat{\mathcal{M}}^d = \check{\mathcal{M}}^d/\Gamma \to \mathcal{A}^0 = \hat{\mathcal{A}}.
\]

All the three Hitchin maps $\chi$, $\cchi$ and $\hchi$ we defined above are proper algebraically completely integrable systems, therefore as explained 
in \S\ref{intro}.\ref{integrablesystem} we should expect the generic fibers to become compact tori, and as we are in the algebraic situation: Abelian varieties. We will see below \S\ref{higgs}.\ref{genericfibers} that this is indeed the case.

\subsection{Cohomology of Higgs moduli spaces} Compared to the moduli spaces of bundles $\N$ we have  less
information on the cohomology of $\M$. Only the case of $n=2$ is understood completely. In this case \cite{HT1} describes the universal generators of $H^*(\M)$ and \cite{HT2} describes the relations
among the generators. Universal generators were found in \cite{markman} for all $n$, but we have not even a conjecture about the ring structure of $H^*(\M)$. 

The Betti numbers of $\M$ are known only when $n=2$ by the work of Hitchin \cite{hitchin} and for $n=3$ by the work of Gothen using a Morse theoretical technique which we will describe below \S\ref{higgs}.\ref{n=2tmsdol} in more detail. There is a conjecture
of the Betti numbers for all $n$ in \cite{hausel-villegas}. But for $n>3$ the Betti numbers are not known, although see \S\ref{conclusion}. 

Crucially for us the analogue of Theorem~\ref{harnar} does not hold. Thus  $\Gamma$ acts non-trivially on $H^*(\cM)$ 
and so $H^*(\cM)\ncong H^*(\hM)$ as was already observed by Hitchin in \cite{hitchin} for $n=2$. One can conjecturally describe the non-trivial part of the $\Gamma$-module $H^*(\cM)$ by using ideas from mirror symmetry. We proceed by detailing these ideas.

\section{Topological mirror symmetry for Higgs bundles}
\label{higgs}

The goal of this section is to
establish the global picture:
\[
\xymatrix{
\check{\MC} \ar[dr]_{\check{\chi}} & & \hat{\MC} \ar[dl]^{\hat{\chi}} \\
& \AC^0
}
\]
where $\check{\chi}$ and $\hat{\chi}$ are the Hitchin maps. The
important point is that the generic fibres are (torsors
for) dual abelian varieties. The reason we want to do this is the following. As we mentioned
in the introduction $\M=\M_\Dol$ is isomorphic with the hyperk\"ahler manifold $\M_\Hit$ in complex structure
$I$. If we change complex structure we also have a moduli space interpretation of $\M_\Hit$ in complex
structre $J$. Namely we can identify \cite[\S 9]{hitchin} the complex manifold $(\M_
\Hit,J)$ with a certain  moduli space $\M_\DR$ of flat connections on $C$. In this complex structure 
\beq{syzmdr}
\xymatrix{
\check{\MC}_{\DR} \ar[dr]_{\check{\chi}} & & \hat{\MC}_{\DR} \ar[dl]^{\hat{\chi}} \\
& \AC^0
}
\eeq
 the fibres of the Hitchin map become dual special Lagrangian tori. This is
the setting proposed by Strominger-Yau-Zaslow \cite{strominger-etal} for mirror symmetry 
as discussed in Subsection~\ref{mirror}. We can thus expect $\cM_\DR$ and $\hM_\DR$ to 
be mirror symmetric. In the physics literature such a mirror symmetry was first suggested by \cite{vafa-etal} in 1994
and more recently by Kapustin-Witten \cite{kapustin-witten} in 2006. Below we will be aiming at checking 
the agreements of certain Hodge numbers of $\cM_\DR$ and $\hM_\DR$ which can be called topological mirror symmetry.

\subsection{Spectral curves}

The simple idea of describing a polynomial by its zeroes leads to the notion of spectral curve of a Higgs bundle $(E,\phi)$ or more generally of its characteristic polynomial $a\in \A$. Recall that it has the form $$a=t^n+a_1 t^{n-1} + \cdots  + a_n,$$ where $a_i\in H^0(K^i)$.
What should be the spectrum of such a polynomial? Look at one point $p \in C$, there we get $\phi_p: E_p \to E_p \otimes K_p$, we expect of an eigenvalue $\nu_p$ of $\phi_p$ to satisfy
	that there exists $0\neq v \in E_p$  with the property $\phi_p(v) = \nu_pv.$
Thus, we need $\nu_p \in K_p$. We do now consider all eigenvalues as a subset of the total space $X$ of the bundle $K \to C$, and want to identify it with the complex points of a scheme. The resulting object will be called the spectral curve corresponding to $a \in \mathcal{A}$ and denoted by $C_a$. The picture is this:
\begin{displaymath}
	\xymatrix{ C_a \ar[rrd]^{\pi_a} &\subset & X \ar[d]^{\pi} \\
			& & C}
\end{displaymath}
To construct the scheme structure on the spectral curve $C_a$, note that there exists a tautological section $\lambda \in H^0(X,\pi_a^*K)$ satisfying $\lambda(x) = x$. We can now pullback the sections $a_i$ to $X$ and obtain a section
\[ s_a = 	\lambda^n + a_1 \lambda^{n-1} + \cdots + a_n \in   H^0(X,\pi_a^*K^n) \]
Clearly $C_a$ equals the zero set of this section, i.e. $C_a = s_a^{-1}(0)$, which comes naturally with a scheme structure.

\subsection{Generic fibres of the Hitchin map}\label{genericfibers}
The fibers of the Hitchin map can be  complicated: reducible even non-reduced. But for generic $a \in \mathcal{A}$ the spectral curve $C_a$ is smooth and the corresponding fiber of the Hitchin map is also smooth and has a nice description. Here we review this description.  For more details
see \cite[\S 8]{Hitchin} and \cite{beauville-etal}.

If $(E,\phi)$ is a Higgs bundle with characteristic polynomial $a$, then the pull-back of $E$ to $C_a$ will have a natural subsheaf $M$ given generically
by the eigenspace $\langle v\rangle\in K_p$ of $\nu_p\in C_a$.  As $C_a$ is non-singular $M$ becomes an invertible sheaf, a line bundle. A more precise definition is to take $$M:={\rm{ker}} \{ \pi^*_a(\phi)-\lambda Id:\pi_a^*(E) \to\pi_a^*(E\otimes K)\}.$$  One can recover $E$ from $M$ by the formula \beq{eigen} E=\pi_{a*}(L),\eeq where $L=M(\Delta)$ and $\Delta$ is the ramification divisor of $\pi_a$ (cf. \cite[Remark 3.7]{beauville-etal}). Then equation \eqref{eigen} can be considered the {\em eigenspace decomposition} of $\phi$ on $E$.    

More generally we  want to get a correspondence between line bundles on $C_a$ and Higgs bundles on $C$ with characteristic polynomial $a$. Starting with a line bundle $L$ on $C_a$ we do at least know that $\pi_{a*}(L)$ is a torsion free sheaf on $C$, but since $C$ is a non-singular curve this means that it is actually a vector bundle of rank $n$, which is the degree of the covering $\pi_a:C_a\to C$. Recall the canonical section $\lambda \in H^0(X,\pi_a^*K)$. It gives us a homomorphism 
\bes L \stackrel{\lambda}{\longrightarrow}   L \otimes \pi_a^*K .
\ees
We can now push this forward to the curve $C$ to obtain
\bes E = \pi_{a*}(L) \stackrel{\pi_{a*}(\lambda)}{\longrightarrow}   \pi_{a*}(L \otimes \pi_a^*K) = \pi_{a*}(L) \otimes K .
\ees
 This way from a line bundle $L$ on $C$ we  get a Higgs bundle $E \to E \otimes K$.

As $C_a$ is integral, $(E,\varphi)$ cannot have any sub Higgs bundle (as its spectral curve would be a one-dimensional subscheme of $C_a$), hence it is automatically stable. Thus if we set $$d':=d+n(n-1)(g-1)$$ then we get a map:
\begin{eqnarray*}
	 &\Pic^{d'}(C_a) \to  \mathcal{M}^{d} \\
			&L \mapsto  (\pi_{a*}(L), \pi_{a*}(\lambda))
\end{eqnarray*}
To see that the Higgs bundle $(\pi_{a*}(L),\pi_{a*}(\lambda))$ has characteristic polynomial $a$ recall that $$\lambda^n + a_{n-1} \lambda^{n-1} + \cdots + a_0 = 0$$ holds on $C_a$. So if we  push-forward this equation to $C$ we obtain $$\pi_{a*}(\lambda)^n + a_{n-1} \pi_{a*}(\lambda)^{n-1} + \cdots + a_0 = 0.$$ Now the Cayley-Hamilton theorem implies the assertion.

\begin{theorem}[Hitchin 1987, Beauville-Narasimhan-Ramanan 1989] \label{hitchinfibres}When $a \in \mathcal{A}_{reg}$ ($\Leftrightarrow$ $C_a$ is non-singular) we have $\chi^{-1}(a) \cong \Pic^{d'}(C_a)$.	
\end{theorem}
We need some modifications for $\SL_n$ and $\PGL_n$. In the $\SL_n$-case we have $a \in \mathcal{A}^0$, we need to find the line bundles $L$, such that $\pi_{a*}(L)$ has the right determinant. Define $\Prym^{d'}(C) \subset \Jac^{d'}(C_a)$ by
\[
	L \in \Prym^{d'}(C_a) \Leftrightarrow \det \pi_{a*}(L)  = \Lambda
\]
It is clear that for all $a \in \mathcal{A}^0_{reg}$, the Hitchin fibre satisfies $\check{\chi}^{-1}(a) \cong \Prym^{d'}(C_a)$.

For $\PGL_n$ we have $\check{\chi}^{-1}(a) \cong \Prym^{d'}(C_a)/\Gamma$. This makes sense since for $L_\gamma \in \Pic(C)[n]$ we do have $$\det (\pi_{a*}(\pi_a^*(L_{\gamma}) \otimes L)) = \det(L_{\gamma} \otimes \pi_{a*}(L)) = L_{\gamma}^n \otimes \det(\pi_{a*}L) = \det(\pi_{a*}L).$$

To summarize, the fibres of the Hitchin map are given:

\begin{itemize}
\item For $\GL_n$: By Thm \ref{hitchinfibres}, for $a \in \AC_{reg}$ 
\[\mathcal{A}_a:=\chi^{-1}(a)\simeq \Jac^{d'}(C_a).\]
\item For $\SL_n$:  following the definitions it is straighforward that for $a \in \mathcal{A}_{reg}^0$ 
\[\check{\mathcal{A}}_a:=\check{\chi}^{-1}(a)\simeq \Prym^{d'}(C_a).\]
\item For $\PGL_n$: There are two ways of thinking of the Hitchin fibre:
\[\hat{\mathcal{A}}_a:=\Prym^{d'}(C)/\Gamma \simeq \Jac^{d'}(C_a)/\Pic^0(C),\]
where $\Pic^0(C)$ acts on $\Jac^{d'}(C_a)$ by tensoring with the pull-back line bundle. 
A short computation shows that the $\Gamma\subset\Pic^0(C)$ action preserves $\Prym^{d'}(C)$. 
\end{itemize}

\subsection{Symmetries of the Hitchin fibration}
\label{symmetries}
We will see in this subsection how natural Abelian varieties act on the regular fibers of the three Hitchin map, giving them a torsor structure. 
Again, we study the $\GL_n, \SL_n, \PGL_n$ cases separately.

 \subsubsection{For $\GL_n$} Fixing $a\in \mathcal{A}_{reg}$, tensor product defines a simply transitive action of $\Pic^0(C_a)$ on $\Jac^d(C_a)$, and therefore $\mathcal{M}_a$ is a torsor for $P_a:=\Pic^0(C_a)$. 
\subsubsection{For $\SL_n$}  Fix $a\in \mathcal{A}_{reg}^0$, we have the (ramified) spectral covering map $\pi:C_a \to C$. 
\begin{definition} The {\it norm map} 
\beq{normap} \Nm_{C_a/C}: \Pic^0(C_a) \to \Pic^0(C)\eeq
is defined by any of the following three equivalent way:
\begin{enumerate}
\item Using divisors. For any divisor $D$ on $C_a$ we have $$\Nm_{C_a/C}(\mathcal{O}(D))=\mathcal{O}(\pi_{a*}D),$$ where $\pi_a: C_a \to C$ is the projection. Then one can show that the norm of a principal divisor will be a principal divisor (using  the norm map between the function fields $\C(C_a)\to\C(C)$ - 
this justifies the name "norm map") - thus inducing a well-defined norm map as in \eqref{normap}.
This definition points out why $\Nm$ is a group homomorphism. 
\item For $L\in \Pic^0(C_a)$ define $$\Nm_{C_a/C}(L)=det(\pi_{a*}(L))\otimes det^{-1}(\pi_{a*}\mathcal{O}_{C_a}).$$
\item Using the fact that $\Pic^0(C), \Pic^0(C_a)$ are Abelian varieties, we can define the norm map as the dual of the pull-back map $$\pi_a^*:\Pic^0(C) \to \Pic^0(C_a),$$ that is
\[\Nm_{C_a/C}=\check{\pi}: \Pic^0(C_a)\simeq\check{\Pic}^0(C_a) \to \check{\Pic}^0(C) \simeq \Pic^0(C).\] Here recall that for an Abelian variety $A$ the dual $\hat{A}$ Abelian variety denotes the moduli space of degree $0$ line bundles on $A$ and that the dual of
the Jacobian of a smooth curve is itself. 
\end{enumerate}
\end{definition}
Let \beq{smoothprym}\Prym^0(C_a):=\ker(\Nm_{C_a/C})\eeq denote the kernel of the norm map, which is an Abelian subvariety of $\Pic^0(C)$. Then $\Prym^0(C_a)$ acts on $\Prym^d(C_a)=\check{\mathcal{M}}_a$, and $\check{\mathcal{M}}_a$ is a torsor for $\check{P}_a:=\Prym^0(C_a)$. 

\subsubsection{For $\PGL_n$} In this case 
\[\hat{\mathcal{M}}_a=\check{\mathcal{M}}_a/\Gamma=\mathcal{M}_a/\Pic^0(C)\] 
is a torsor for $\hat{P}_a:=\check{P}_a/\Gamma=P_a/\Pic^0(C)$.

To complete the SYZ picture \eqref{syzmdr} we show that $\check{P}_a$ and $\hat{P}_a$ are dual abelian varieties. 

\subsection{Duality of the Hitchin fibres}

Take the short exact sequence
\[\begin{array}{ccccccccc}
0 & \rightarrow & \Prym^0(C_a) & \hookrightarrow & \Pic^0(C_a) & \stackrel{\Nm_{C_a/C}}{\twoheadrightarrow}& \Pic^0(C) & \rightarrow & 0
\end{array}\]
and dualize. Since $\Pic^0(C)$ and $\Pic^0(C_a)$ are isomorphic to their duals, we get
 \[\begin{array}{ccccccccc}
0 & \leftarrow & \check{\Prym}^0(C_a) & \twoheadleftarrow & \Pic^0(C_a)\stackrel{\pi_a^*}{\leftarrow} & \Pic^0(C) & \leftarrow & 0
\end{array},\]
and therefore 
\[\check{\check{P}}_a=\Pic^0(C_a)/\Pic^0(C)=\hat{P}_a,\]
that is $\check{P}_a$ and $\hat{P}_a$ are duals. (See \cite[Lemma (2.3)]{hausel-thaddeus} for more details.)

This is the first reflection of mirror symmetry. To summarize, we can state
\begin{theorem}[\cite{hausel-thaddeus}]
For a regular $a \in \mathcal{A}_{reg}^0$ the fibers $\check{\mathcal{M}}_a$ and $\hat{\mathcal{M}}_a$ are torsors for dual Abelian varieties $\check{P}_a$ and $\hat{P}_a$, respectively. 
\end{theorem}

We can state this theorem more precisely using the language of gerbes. To that end here is a short summary.

\subsection{Gerbes on $\check{\mathcal{M}}$ and $\hat{\mathcal{M}}$} \label{gerbes}
Here we sketch a quick definition of gerbes for more details see \cite{hausel-thaddeus},\cite{thaddeus} and \cite{donagi-gaitsgory}.
Let $A$ be a sheaf of Abelian groups on a variety $X$. The typical examples are $\mathcal{O}_X^\times$, and the constant sheaves $\underline{\mu_n}, \underline{\U(1)}$, where $\mu_n$ is the group of $m$th roots of unity. Note that $\underline{\mu_n} \subset \mathcal{O}_X^\times$ and $\underline{\mu_n} \subset \U(1)$.
\begin{definition}
An $A$-torsor is a sheaf $F$ of sets on $X$ together with an action of $A$, such that $F$ is locally isomorphic with $A$.  In particular, when nonempty, $\Gamma(U,F)$  is a torsor for $\Gamma(U,A)$ for all open $U\subset X$.
\end{definition}
Examples: \begin{itemize} \item
 $\mathcal{O}_X^\times$-torsor =  line bundle  \item
 $\underline{\U(1)}$-torsor = flat unitary line bundle \item
 $\underline{\mu_n}$-torsor = $\mu_n$-Galois cover 
\end{itemize}
Note that the natural tensor category structure on the category of torsors $Tors_A(U)$ endows it with a group-like structure. Moreover, the 
automorphism of an $A$-torsor is an element of $\Gamma(A)$. 

\begin{definition}
An $A$-gerbe $B$ is a sheaf of categories (which is roughly what one would think) so that locally $B|_U$ becomes the analogue of  a torsor over $Tors_A(U)$. 
\end{definition}

Let $(\mathbb{E},{\Phi})$ be a universal Higgs-bundle on
$\check{\mathcal{M}}\times C$, where $${\Phi}\in
H^0(\End_0\mathbb{E}\otimes \pi_a^*(K_C)),$$ and $\mathbb{E}_c=\mathbb{E}|_{\check{\mathcal{M}}\times \{c\}}$ be the
fiber over $c \in C$. Such a universal bundle exists  by our running assumption $(d,n)=1$.  Then 
\begin{equation} \label{generator}  c_1(\mathbb{E}_c)\in H^2(\cM,\mathbb{Z})\simeq \mathbb{Z} {\mbox{ is a generator modulo $n$}}.  \end{equation} Note that $\mathbb{E}$ is not unique: it can be tensored by $L\in \Pic(\cM)$, but this property always holds. 

Let $\mathbb{P}\mathbb{E}_c \to \check{\mathcal{M}}$ be the corresponding $\PGL_n$-bundle. Let $\check{B}$ be the $\mu_n$-gerbe of liftings of $\mathbb{P}\mathbb{E}_c$ as an $\SL_n$-bundle. Because for every lifting as a $\GL_n$-bundle \eqref{generator} holds, there is no global lifting as an $\SL_n$-bundle and so $\check{B}$ is not a trivial gerbe. 
But for $a \in \mathcal{A}_{reg}$ it turns out that
$c_1(\mathbb{P}\mathbb{E}_c)|_{{\check{\mathcal{M}}_a}}$ is $0$ mod $n$, and so $\mathbb{P}\mathbb{E}_c|_{\check{\mathcal{M}}_a}$ can be lifted as an $\SL_n$-bundle, therefore $\check{B}|_{\check{\mathcal{M}}_a}$ is a trivial gerbe. 

Finally we note that the action of $\Gamma$ on $\check{\M}$ can be lifted to $\check{B}$ to get an orbifold gerbe $\hat{B}$ on $\hat{\M}$. As $\check{B}$
and $\hat{B}$ are $\mu_m\subset \U(1)$-gerbes we can consider the induced $\U(1)$-gerbes as well. 

\begin{theorem}[\cite{hausel-thaddeus}] \label{fullsyz} One can identify the set of 
	trivializations
\[Triv^{\U(1)}(\check{B}^e|_{\check{\mathcal{M}}_a^d})\simeq \hat{\mathcal{M}}_a^e\]
as $\check{P}_a$-torsors. 
Similarly, 
\[
Triv^{\U(1)}(\hat{B}^d|_{\hat{\mathcal{M}}_a^e}) \simeq
\check{\mathcal{M}}_a^d.
\]
\end{theorem}

\begin{remark} 
The analogue of this result for arbitrary pair of Langlands dual groups was handled by Donagi--Pantev in \cite{donagi-pantev1}. There they also implement
the fiberwise Fourier-Mukai transform over the regular locus which will be discussed in \S\ref{solving}.\ref{shadow}. The case of
${\rm G}_2$ was considered by Hitchin in \cite{hitchin4} in detail. 

This Theorem~\ref{fullsyz}  can be interpreted as the twisted version of the Strominger-Yau-Zaslow proposal suggested in \cite{hitchin-twistSYZ}. Thus
we have established the picture

\bes
\xymatrix{
\check{\MC}^d_\DR \ar[dr]_{\check{\chi}} & & \hat{\MC}^e_\DR \ar[dl]^{\hat{\chi}} 
\\
& \AC^0}
\ees 
with $\check{\chi}$ and $\hat{\chi}$ two special Lagrangian fibrations
with dual tori as generic fibers according to Theorem~\ref{fullsyz}. Thus 
the pair $(\cM_\DR^d,\check{B}^e)$ and $(\cM_\DR^e,\hat{B}^d)$ can considered
mirror symmetric in the twisted SYZ sense of \cite{hitchin-twistSYZ} Hitchin. 

\end{remark}

\subsection{The stringy Serre polynomial of an orbifold}

As our spaces satisfy a suitable version of the Strominger-Yau-Zaslow mirror symmetry  proposal 
we need a definition of Hodge numbers of non-projective varieties to be able 
to formulate  topological mirror symmetry for our  mirror pairs. For more details on these Hodge numbers see \cite{deligne2, batyrev-dais,hausel-thaddeus, hausel-villegas, peters-steenbrink}.

\begin{definition}\label{epoly}
Let $X$ be a complex algebraic variety. The Serre polynomial (virtual Hodge polynomial) is defined as
\[E(X;u,v)=\sum_{p,q,i}(-1)^ih^{p,q}(Gr_{p+q}^WH_c^i(X))u^pv^q\]
where $Gr_{p+q}^WH_c^i(X)$ is the $p+q$th graded piece of the weight filtration of $H_c^i(X)$. This has a Hodge structure of weight $p+q$, and $h^{p,q}$ is the corresponding Hodge number. This so-called mixed Hodge structure was constructed by Deligne \cite{deligne2}.
\end{definition}
\begin{remark}
We say that the mixed Hodge structure is  pure if $h^{p,q}\neq 0$ implies that $p+q=i$. For example, smooth projective varieties have pure MHS on their cohomology. One can also show \cite[Theorems 2.1,2.2]{hausel-mln} that the spaces $\M^d_\Dol$ and $\M^d_\DR$ also have pure MHS on their cohomology.  
\end{remark} The definition of the virtual Hodge polynomial is the first step of formulating the topological mirror test. In order to obtain the correct Hodge numbers of the orbifold  $\hM$. We need to define stringy Hodge numbers twisted by a gerbe. 
Let $X$ be smooth, and assume that a finite group $\Gamma$ acts on $X$. Then we can define the \textsl{stringy Serre polynomial} of the orbifold $X/\Gamma$ as follows:
\begin{equation}\label{stringyorb}
E_\st(X/\Gamma),u,v)=\sum_{[\gamma]\in [\Gamma]}E(X_\gamma/C_\gamma;u,v)(uv)^{F(\gamma)}, \end{equation}
\vskip-.4cm where 
\begin{itemize}
 \item $[\gamma]$ is a conjugacy class of $\Gamma$ \item $X_\gamma$ is the fixed point set, $C_\gamma$ is the centralizer of $\gamma$ in $\Gamma$, acting on $X_\gamma$. 
\item $F(\gamma)$ is the Fermionic shift, defined as $F(\gamma)=\sum w_i$, where $\gamma$ acts on $TX|_{X_\gamma}$ with eigenvalues $e^{2\pi iw_i}$, $w_i \in [0,1)$. 
\end{itemize}
\begin{remark} These orbifold Hodge numbers can be considered the cohomological shadow of the stringy derived category of coherent sheaves on the variety $X/\Gamma$, or more topologically the stringy $K$-theory of $X/\Gamma$. Both can be simply defined by  considering $\Gamma$-equivariant coherent sheaves  or $\Gamma$-equivariant complex vector bundles on $X$. See \cite{jarvis-etal} for more details. 
\end{remark}
Another important property of $E_{\rm st}(X/\Gamma;u,v)$ is the following theorem:

\begin{theorem}[\cite{kontsevich2,batyrev-dais}] If $\pi:Y\to X/\Gamma$ is a crepant resolution, i.e. $\pi^*\omega_{X/\Gamma}=\omega_Y$, then $E_{\rm st} (X/\Gamma;u,v)=E(Y;u,v)$.\end{theorem} 

\begin{remark} This is a nice way to see what the stringy Hodge numbers mean, in terms of the Hodge numbers of (any) crepant resolution. However our orbifolds never have crepant resolutions, because the generic singular points are infinitesimally modelled by quotients of symplectic representations of $\Z_n$ which are not generated by symplectic reflections 
as in \cite[Theorem 1.2]{verbitsky}.  So we will be stuck with the definition in \eqref{stringyorb}.
\end{remark}

Finally let $B$ be a $\Gamma$-equivariant $\U(1)$-gerbe on $X$.  Then more generally, we define
\begin{equation}
E_\st^B(X/\Gamma;u,v)=\sum_{[\gamma]\in [\Gamma]}E(X_\gamma/C_\gamma,L_{B,\gamma};u,v)(uv)^{F(\gamma)},
\end{equation}
where $L_{B,\gamma}$ is the local system on $X_{\gamma}$ given by $B$. 

These $B$-twisted stringy Hodge numbers can again be considered as the cohomological shadow of a certain derived category of  $B$-twisted $\Gamma$-equvariant sheaves on $X$ or the corresponding $B$-twisted equivariant $K$-theory.

\subsection{Topological mirror test}

We can now formulate our original topological mirror symmetry

\begin{conjecture}[\cite{hausel-thaddeus}]
\label{drtms}
For $(d,n)=(e,n)=1$
\[E(\check{\mathcal{M}}_\DR^d;u,v)=E_\st^{\hat{B}^d}(\hat{\mathcal{M}}_\DR^e;u,v),\] 
where $\hat{B}^d$ is the $\Gamma$-equivariant gerbe on $\check{\M}^e$ appearing in Theorem~\ref{fullsyz}.
\end{conjecture}

\begin{remark}\label{handwaving} One can object that this conjecture does not feature the change in the indices of the Hodge numbers as in the original topological mirror test \eqref{otms} between projective Calabi-Yau manifolds. Our, rather hand-waving, argument for this in \cite{hausel-thaddeus} was that 
for {\em compact} hyperk\"ahler manifolds \eqref{otms} is equivalent with \eqref{ctms} that is agreement of Hodge numbers. This was only what we could offer to explain the agreement of Hodge numbers in Conjecture~\ref{drtms},  even though our examples are not compact, and their Hodge numbers do not possess any non-trivial  symmetries. As will be explained after Conjecture~\ref{tms}  some recent developments lead to a solution of  this problem. 
\end{remark}

As   $\M_\DR$ can be algebraically deformed to  $\M_\Dol$, inside nice compactifications,  it is not  surprising
that we could prove that their mixed Hodge structure is isomorphic:
\begin{theorem}[{\cite[Theorem 6.2,6.3]{hausel-thaddeus}} ]\label{dr=dol} For $(d,n)=
(e,n)=1$ $$E(\check{\mathcal{M}}_\DR^d;u,v)=E(\check{\mathcal{M}}_\Dol^d;u,v)$$ and similarly $$E_\st^{\hat{B}^d}(\hat{\mathcal{M}}_\DR^e;u,v)=E_\st^{\hat{B}^d}(\hat{\mathcal{M}}_\Dol^e;u,v).$$
\end{theorem}

Thus we have an equivalent form of Conjecture~\ref{drtms}, the so-called Dolbeault version of topological mirror symmetry:
\begin{conjecture}[\cite{hausel-thaddeus}]
\label{doltms}
For $(d,n)=(e,n)=1$
\[E(\check{\mathcal{M}}_\Dol^d;u,v)=E_\st^{\hat{B}^d}(\hat{\mathcal{M}}_\Dol^e;u,v),\] 
where $\hat{B}^d$ is the $\Gamma$-equivariant gerbe on $\check{\M}^e$ appearing in Theorem~\ref{fullsyz}.
\end{conjecture}
As we will sketch in \S\ref{solving}.\ref{shadow} this conjecture could be interpreted as a  cohomological shadow of some equivalence of derived categories of sheaves on the Hitchin systems
which arise in the Geometric Langlands programme.

This latter conjecture we were able to prove in the following cases:

\begin{theorem}[\cite{hausel-thaddeus}] Conjecture~\ref{doltms} and so Conjecture~\ref{drtms} are valid for
$n=2$ and $n=3$.
\end{theorem}

In \S\ref{higgs}.\ref{n=2tmsdol} below we will sketch the computation in the $n=2$ case. 

Now  we will unravel the meaning of Conjecture~\ref{doltms}.
Recall that $\Gamma$ acts on $\cM^d$ and hence $\Gamma$ also acts
on $H_c^*(\cM^d)$. We get a decomposition
\[
H_c^*(\cM^d) = \bigoplus_{\kappa \in \hat{\Gamma}} H^*_c(\cM^d)_{\kappa}
\]
which is compatible with the mixed Hodge structure. Therefore we can
write

\beq{edecomp}E( \cM^d; u, v) = \sum_{\kappa \in \hat{\Gamma}}
E_{\kappa}(\cM^d; u,v),\eeq where $$E_{\kappa}(\cM^d; u,v):=\sum_{p,q,i}(-1)^ih^{p,q}(Gr_{p+q}^WH_c^i(\cM^d)_\kappa)u^pv^q.$$

We can also expand the RHS of Conjecture~\ref{doltms}. By the definition of the  stringy Serre polynomial, and because $\Gamma$ is commutative,  
we have
\[
E_\st^{\hat{B}^d}( \hM^e; u,v) = 
 \sum_{\gamma \in \Gamma}
E(\cM^e_{\gamma}/ \Gamma; L_{\hat{B}^d,\gamma}, u, v)(uv)^{F(\gamma)}
\]

Thus the unravelled Conjecture~\ref{doltms} takes the form:

\beq{unravel}\sum_{\kappa \in \hat{\Gamma}}
E_{\kappa}(\cM^d; u,v)=\sum_{\gamma \in \Gamma}
E(\cM^e_{\gamma}/ \Gamma; L_{\hat{B}^d,\gamma}, u, v)(uv)^{F(\gamma)}
\eeq

We note that there are the same number of terms in \eqref{unravel}; and in fact there is  a canonical way to identify them.  Note that $\Gamma$ is canonically isomorphic to $H^1(C, \Z_n)$, where
$C$ is our underlying curve. It follows that Poincar\'e duality gives
us a canonical pairing
\[
w : \Gamma \times \Gamma \to H^2(C, \Z_n) = \Z_n,
\] the so-called Weil pairing.
This allows us to identify $w : \Gamma \to \hat{\Gamma}$.
This identification  leads to the \emph{refined topological mirror symmetry test}:

\begin{conjecture}\label{refinedtms}For $\kappa\in \hat{\Gamma}$ we have \[
E_{\kappa}( \cM^d; u, v) = E( \cM^e_{\gamma}/ \Gamma,
L_{\hat{B}^d, \gamma} ; u, v)(uv)^{F(\gamma)}
\] 
where $\gamma = w(\kappa)$. \end{conjecture} This again holds for $n=2,3$ the case of $n=2$ will be discussed in the  next section.   As we will see in
\S\ref{solving}.\ref{TMStoFL} this refined conjecture is closely related to Ng\^{o}'s geometric approach to the fundamental lemma.

Here we point out the case of the trivial character $\kappa=1$. In that case the refined topological mirror symmetry test is:

\begin{conjecture} \label{e=d} When $(d,n)=(e,n)=1$ 
$$E_{1}(\cM^d;u,v)=E(\cM^e/\Gamma;u,v)$$ or equivalently 
$$E(\hM^d;u,v)=E(\hM^e;u,v).$$
\end{conjecture}
\begin{remark} This conjecture is in sharp contrast to the dependence on degree of $H^*(\hN^d)$. In fact it is known \cite{zagier} that already when $n=5$ we have $$E(\hN^1;u,v)\neq E(\hN^3;u,v).$$ 
\end{remark}
\begin{remark}
As Theorem~\ref{fullsyz} holds for all $d$ and $e$ not necessarily coprime to $n$ it is concievable that something like Conjecture~\ref{e=d} should hold even when $d$ (and/or $e$)  is not coprime to $n$. It is however unclear what cohomology theory we should 
calculate on the singular moduli space $\cM^d$, for $(d,n)\neq 1$, to produce the agreement in 
Conjecture~\ref{e=d}. For $n=2$ Batyrev's extension \cite{batyrev} of stringy cohomology of $\cM^0$ was calculated in \cite{kim},  while stacky cohomology of $\cM^0$ was calculated in \cite{wilkin-etal}. In either case 
the corresponding generating functions are not polynomials; thus Conjecture~\ref{e=d} for $n=2, d=0,e=1$, as it stands, cannot hold for them. 
\end{remark}



\subsection{Topological mirror symmetry for $n = 2$ }\label{n=2tmsdol}

 Consider the circle action of $\C^\times$ on the Higgs moduli space by rescaling the Higgs field. That
is, $\lambda \cdot (E, \phi) \mapsto (E, \lambda \cdot \phi)$. We can
study the corresponding Morse stratification and by Morse theory we obtain the
decomposition
\beq{morse}
H^*( \cM) = \bigoplus_{F_i\subset \cM^{\C^\times}} H^{* + \mu_i}(F_i).
\eeq
Here the sum is over the connected components of the fixed point set
$\cM^{\C^\times}$, and $\mu_i$ denotes the index of $F_i$ with
respect to the $\C^\times$-action. Note that \eqref{morse} is a
decomposition as $\Gamma$-modules.

The components $F_i$ have been described for $n = 2$ by Hitchin \cite{hitchin},
and by Gothen \cite{gothen} for $n =3$. The case $n = 4$ seems quite hard; but
for recent progress see \S\ref{conclusion}.

One obvious fixed point locus is $F_0\cong \cN$, consisting of stable bundles
with zero Higgs field. However this component doesn't contribute to the variant
part as the $\Gamma$-action on $H^*(\cN)$ is trivial by Theorem~\ref{harnar}.

From now on in this section we assume $n=2$ and $d=1$. Then  the other components can be labelled by $i= 1, \dots, g-1$ and consist of
isomorphism classes of Higgs bundles of the form
\[
F_i = \{ (E, \phi) \; | \; E \cong L_1 \oplus L_2, \deg(L_1)=i, \phi = \left
 ( \begin{array}{cc} 0 & 0 \\  \varphi & 0 \end{array} \right ),
0\neq \varphi \in H^0( L_1^{-1}L_2K)
\}.
\]
Now stability forces $\deg L_2 = 1-i$ where $i >0$, because $L_2$ is a
$\phi$-invariant subbundle. We can associate to $(E, \phi) \in F_i$
the divisor of $\varphi$ in $S^{2g-2i-1}(C)$ yielding a $2^{2g}:1$ covering
\[
F_i \to S^{2g-2i-1}(C)
.\] This covering is given by the free action of $\Gamma$ on $F_i$.

\begin{theorem}[\cite{hitchin}(7.13)]
The $\Gamma$-action on $H^*(F_i)$ is only non-trivial in the middle
degree $2g -2i -1$. We have
\begin{gather*}
\dim H^{2g - 2i - 1}_{var}(F_i) = (2^g - 1) {{2g - 2} \choose {2g - 2i
 - 1}}
\end{gather*}
Moreover, if $\kappa \in \hat{\Gamma}^*$ then
\begin{gather*}
\dim H^{2g - 2i - 1}_{\kappa}(F_i) = {{2g - 2} \choose {2g - 2i
 - 1}}.
\end{gather*}
\end{theorem}

We now consider the stringy side. Recall $\hat{\MC} = \check{\MC}/\Gamma$ and let
$\gamma \in \Gamma^*=\Gamma\setminus \{0\}$. Then $\gamma$ leads to a connected covering
\[
\pi_\gamma:C_{\gamma} \to C
\]
with Galois group $\Z _2$. Consider the
commutative diagram
\[
\xymatrix@C=0.8cm{
 T^* \Jac^d(C_{\gamma})
\ar[r]^(0.6){(\pi_\gamma)_*} \ar[d]^{\cong} & \MC^d \ar[d]^{\det} & \supset \check{\MC}^d =
\det^{-1}(\Lambda,0)  \\
T^* \Jac^d(C_{\gamma}) \ar[r]^{\Nm_{C_{\gamma}/C}} & T^* \Jac^d(C) &
\ni (\Lambda,0)
}
\]
From this diagram \cite[Corollary 7.3]{hausel-thaddeus} we have that 
$\check{\MC_{\gamma}}$ is a torsor for  $T^* \Prym(C_{\gamma}/C)^0$; where 
 $\Prym(C_{\gamma}/C)^0:=(\Nm_{C_{\gamma}/C}^{-1}(\calO_C))^0$ is the connected component of the Prym variety. 
We can then calculate that
\[
\dim H^{2g-2i+1}(\check{\MC}_{\gamma} /\Gamma, L_{\hat{B}^d, \gamma}) = { 2g -2
 \choose 2g - 2i - 1},
\] when $i=1,\dots,g-1$
and is zero otherwise. Note that the presence of the gerbe means that
we see only the odd degrees. It follows that in this case we   indeed have the refined topological mirror symmetry

\beq{refinedn=2}
E_{\kappa}(\check{\MC}; u, v ) = E( \check{\MC}_{\gamma}/\Gamma , L_{\hat{B},
 \gamma}; u, v)(uv)^{F(\gamma)}
\eeq
when $\kappa = w^{-1}(\gamma)$. \begin{remark} What is surprising
about the agreement in \eqref{refinedn=2} is that the left hand side comes from the fixed point components $F_i$ for $i>0$, because $F_0=\N$ has no variant cohomology by Theorem~\ref{harnar}. While the right hand side comes  solely $F_0$, because the $\Gamma$ action is free on $F_i$ for $i>0$. Thus in order for
\eqref{refinedn=2} to hold there is a remarkable agreement of cohomological data coming on one hand from the symmetric products $\{F_i\}_{i>0}$ of the curve and 
on the other hand from the moduli space of stable bundles $\cN_\gamma\subset \cN=F_0$.
\end{remark}

One may give a similar proof for \eqref{refinedn=2}, when $n = 3$ using Gothen's work \cite{gothen}. The proof along these lines  for higher $n$ remains incomplete. 

We now introduce a more successful arithmetic technique to prove a different version of topological mirror symmetry. 

\section{Topological mirror symmetry for character varieties}
\label{betti}
As mentioned above in another complex structure $\M_\Dol$ can be identified with  $\M_\DR$ a certain moduli space of twisted flat $\GL_n$-connections on $C$.
 A point in $\M_\DR$ represents a certain twisted flat connection on a rank $n$
bundle. One can take its monodromy yielding a twisted representation of the fundamental group
of $C$. This leads to the character variety of the Riemann surface $C$. 

\begin{definition} The \emph{character variety for $\GL_n$} is defined as the affine
	GIT quotient:
\[
\M_\B^d := \{ (A_1, B_1, \dots, A_g, B_g) \in \GL_n(\C) \; | \; 
[A_1, B_1]\dots [A_g, B_g] = \zeta_n^d I_n \} /\!/ \GL_n
\]
The \emph{character variety for $\SL_n$} is the space
\[
\cM_\B^d := \{ (A_1, B_1, \dots, A_g, B_g) \in \SL_n(\C) \; | \; 
[A_1, B_1]\dots [A_g, B_g] = \zeta_n^d I_n \} /\!/ \GL_n
\]
where the action is always by simultaneous conjugation on all factors, and
$\zeta_n = e^{2\pi i/n}$.
The \emph{character variety for $\PGL_n$} is defined as
\[
{\hM^d_\B} := {\cM^d_\B}/\Gamma = \M_\B^d / (\C^\times)^{2g}.
\] Here $\Gamma := \Z_n^{2g} \subset
(\C^\times)^{2g}$ and $\Gamma\subset(\C^\times)^{2g}$ acts on $\cM^d_\B\subset\M^d_\B$
by multiplying the matrices $A_i,B_i$ by scalars. 
\end{definition}

\begin{remark} For the $\PGL_n$-character variety we could have started to consider  
 $$\hM_\B:=\{ (A_1, B_1, \dots, A_g, B_g) \in \PGL_n(\C) \; | \; 
[A_1, B_1]\dots [A_g, B_g] = I_n \} /\!/ \PGL_n.$$ This variety has $n$ components, depending
on which order $n$ central element in $\GL_n$ will agree with the product of the commutators of the $\GL_n$-representatives
of the $\PGL_n$ elements $A_i,B_i$. For $\hM^d_\B$ we picked the component corresponding to the central element $\zeta^d_n I_n$.  
\end{remark}

One can prove \cite{hausel-villegas} that if $d$ and $n$ are coprime, $\M_\B^d$ and $\cM_\B^d$ are non-singular,
and so $\hM^d_\B$ is an orbifold.

Taking the monodromy will give a complex analytical isomorphism $\M_{\DR}^d\cong\M^d_{\B}$ the
so-called Riemann-Hilbert correspondence. We have the following more general \emph{non-abelian Hodge theorem}
\begin{theorem}[\cite{simpson,corlette,donaldson}] There are canonical diffeomorphisms  in all our cases $\GL_n$, $\SL_n$ and $\PGL_n$: $$\hat{\cM}^d_\Dol{\cong}\hat{\cM}^d_\DR {\cong} \hat{\cM}^d_\B$$
\end{theorem}

Because $\M_\DR$ and $\M_\B$ are analytically isomorphic it follows that the twisted SYZ mirror symmetry proposal is satisfied by the pair $\cM^d_\B$ and ${\hM}^e_\B$ as well. 
We may thus formulate the Betti-version of the topological mirror symmetry
conjecture \cite{hausel-mln}:

\begin{conjecture}\label{btms} For $(d,n)=(e,n)=1$
$$E(\cM_\B^d;u,v)=E_\st^{\hat{B}^d}(\hM_\B^e;u,v),$$
where $\hat{B}$ is the   $\Gamma$-equivariant gerbe on ${\cM_\B}^e$ analogous to the one in \S\ref{higgs}.\ref{gerbes}.
\end{conjecture}
Note that the mixed Hodge structures
on $H^*(\M_\B)\cong H^*(\M_\DR)$ are different, in particular unlike the MHS on  $H^*(\M_\B)$, the MHS on $H^*(\M_\DR)$ is pure. Thus Conjecture~\ref{btms} 
is  different  from Conjecture~\ref{drtms} and the equivalent Conjecture~\ref{doltms}. 

We will also need the  refined version:
\begin{conjecture}\label{bettirefined} For $\kappa\in \hat{\Gamma}$ we have \[
E_{\kappa}( \cM^d_\B; u, v) = E( \cM^e_{\gamma}/ \Gamma,
L_{\hat{B}^d, \gamma} ; u, v)(uv)^{F(\gamma)}
\] 
where $\gamma = w(\kappa)$.
\end{conjecture}
Interestingly,  for $(n,e)=(n,d)=1$  the character varieties $\hM^d_\B$ and $\hM^e_\B$ are Galois conjugate via an automorphism  of the complex numbers sending $\zeta_n^d$ to $\zeta_n^e$. This Galois conjugation induces an isomorphism \beq{galois} H^*(\hM^d_\B)\cong H^*(\hM^e_\B) \eeq
preserving mixed Hodge structures. Therefore the $\kappa=1$ case of refined topological mirror symmetry for character varieties follows:

\begin{theorem} When $(d,n)=(e,n)=1$ 
$$E_{1}(\cM_\B^d/;u,v)=E(\cM_\B^e/\Gamma;u,v)$$ or equivalently 
$$E(\hM_\B^d;u,v)=E(\hM_\B^e;u,v).$$\label{independentd}
\end{theorem}

\begin{remark} In fact one can prove that the universal generators of $H^*(\hM^d_\B)$ 
are mapped to the corresponding ones of $H^*(\hM^e_\B)$, and consequently one can prove that the Galois conjugation \eqref{galois} preserves even the mixed Hodge structure on $H^*(\hM_\Dol^d)\cong H^*(\hM_\Dol^e)$; which implies Conjecture~\ref{e=d}. Thus in the original topological mirror symmetry Conjecture~\ref{doltms} we see that  the first terms corresponding to the trivial elements in $\Gamma$ and $\hat{\Gamma}$ at least agree; so we can concentrate on proving Conjecture~\ref{refinedtms} for non-trivial characters $\kappa$. 
This is the first non-trivial application of considering character varieties. 
\end{remark}

\begin{remark} One application in \cite[Theorem 3.3.2]{harder-narasimhan} was to show that for $(d,n)=(e,n)=1$ the Betti numbers of  $\hN^d$ and $\hN^e$ agree only when $d+e$ or $d-e$ is divisible by $n$. This is markedly different behaviour from Theorem~\ref{independentd} which shows that the Betti numbers of $\hM^d$ and $\hM^e$ agree as long as $(d,n)=(e,n)=1$.
 \end{remark}

In the next section we also offer an arithmetic technique which can be used efficiently to check Conjecture~\ref{btms} in  all cases.

\subsection{An arithmetic technique to calculate Serre polynomials}

Recall the definition of the Serre polynomial of a complex variety $X$:
\[E(X;u,v)=\sum_{i,p,q}(-1)^ih^{p,q}(Gr_{p+q}^{W}H_c^i(X))u^pv^q\]
where $W_0 \subseteq W_1 \subseteq \ldots \subseteq W_i \subseteq \ldots \subseteq W_{2k}=H^k(X):=H^k(X;\Q)$ is the weight filtration.

By \cite[Corollary 4.1.11]{hausel-villegas} $\M_B$ has a Hodge-Tate type MHS, that is, $h^{p,q}\neq 0$ unless $p=q$ in its MHS. In this case the Serre polynomial is a polynomial of $uv$, i.e
\begin{equation}\label{tate}
E(X;u,v)=E(X,uv):=\sum_{i,k}(-1)^i \dim(Gr_k^W H_c^i(X))(uv)^k,
\end{equation}
but the MHS is not pure, i.e there is $k\neq i$ such that $h^{(k/2,k/2)}\neq 0$. 

Roughly speaking a variety $X$ can be defined over the integers if one 
can arrange the defining equations to have integer coefficients. Then one
can consider those equations in finite fields $\F_q$, and can count these solutions to
get the number
$ |X(\mathbb{F}_q)|$. We say that such a variety $X$ is polynomial count if 
$ |X(\mathbb{F}_q)|$
is polynomial in $q$. One can define polynomial count varieties even if they can only be defined over more general finitely generated rings than $\Z$. For more technical details and precise statements see \cite{hausel-villegas}. Here we have 
\begin{theorem}[{\cite[Appendix by Katz]{hausel-villegas}}]\label{katz}
For a polynomial count variety $X$ 
\[E(X/\mathbb{C},q)=|X(\mathbb{F}_q)|.\]
\end{theorem}

\begin{example}
Define $\mathbb{C}^*=\mathbb{C}\backslash \left\{0\right\}$ over $\mathbb{Z}$ as the subscheme $\left\{xy=1\right\}$ of $\mathbb{A}^2$. Then 
\[E(\mathbb{C}^*;q)=|\mathbb{F}_q^*|=q-1.\]
Since $H_c^2(\mathbb{C}^*)$ has weight $2$ and $H_c^1(\mathbb{C}^*)$ has weight $0$, substitution to \eqref{tate} gives indeed $q-1$ for $E(\mathbb{C}^*;q)$. 
\end{example}

\subsection{Arithmetic harmonic analysis on $\MC_B$}

By Fourier transform on a finite group $G$ one gets the following Frobenius-type formula:
\[ \left| \left\{a_1,b_1,\ldots ,a_g,b_g \in G|\prod [a_i,b_i]=z\right\}\right|=\sum_{\chi \in Irr(G)}\frac{|G|^{2g-1}}{\chi(1)^{2g-1}}\chi(z).\]
Therefore assuming that $\zeta_n \in \mathbb{F}_q^*$, i.e $n|q-1$, we get
\begin{equation}\label{nofpoints}
| \M_B^d(\mathbb{F}_q) |=(q-1)\sum_{\chi \in Irr(\GL_n(\mathbb{F}_q))}\frac{|\GL_n(\mathbb{F}_q)|^{2g-2}}{\chi(1)^{2g-2}}\cdot \frac{\chi(\zeta_n^d\cdot I)}{\chi(1)}.
\end{equation}

Irreducible characters of $\GL_n(\mathbb{F}_q)$ have a combinatorial description by
Green \cite{green} from 1955. Consequently  $| \M_B^d(\mathbb{F}_q) |$ can be calculated explicitly \cite{hausel-villegas}. It   turns out to be a polynomial, so Katz's Theorem~\ref{katz} applies and \eqref{nofpoints} gives the Serre polynomial. 

The same Frobenius-type formula is valid in the $\SL_n$-case:
\begin{equation}\label{nofpointssl}
| \check{\MC}_B^d(\mathbb{F}_q)|=\sum_{\chi \in Irr(\GL_n(\mathbb{F}_q))}\frac{|\SL_n(\mathbb{F}_q)|^{2g-2}}{\chi(1)^{2g-2}}\cdot \frac{\chi(\zeta_n^d\cdot I)}{\chi(1)}, 
\end{equation}
Here the character table of $\SL_n(\F_q)$ is trickier. After much
work of Lusztig, the character table of
$\mathrm{Irr}(\SL_n(\mathbb{F}_q))$ has only been completed by
Bonnaf\'e \cite{bonnafe} and Shoji \cite{shoji} in 2006.  However for $\chi \in
\mathrm{Irr}(\GL_n(\mathbb{F}_q))$ the splitting
\[\chi|_{\SL_n(\mathbb{F}_q)}=\sum \chi_i\] into irreducible characters $\chi_i$ of $\SL_n(\F_q)$ 
is evenly spread out on $\zeta_n^d \cdot I$; meaning that $\chi_i(\zeta_n^d \cdot I)=\chi_j(\zeta_n^d \cdot I)$. This way the evaluation of \eqref{nofpointssl}   is possible and was done by  Mereb \cite{mereb} in 2010. He too obtained 
a polynomial 
for $|\check{\MC}_B^d(\mathbb{F}_q)|$, and by Katz's theorem Theorem~\ref{katz}
this gives a formula for $E(\check{\MC}_B^d(\mathbb{F}_q);q)$. 
With similar techniques one can also evaluate the $\kappa$-components $E_\kappa(\cM^d_\B;u,v)$ in the LHS of Conjecture~\ref{bettirefined}. 

In order to check the refined topological mirror symmetry Conjecture~\ref{bettirefined} we need also to determine $E(\hat{\MC}_{B,\gamma}^e,L_{\hat{B}^d,\gamma};q)$. It is simple to do when $n$ is a prime; leading to a proof of Conjecture~\ref{bettirefined} and so to Conjecture~\ref{btms} in this case. For composite $n$'s an ongoing work \cite{hausel-etal}  evaluates these by similar (twisted) arithmetic techniques. This seems to match with Mereb's result, which is expected to give the  proof of Conjecture~\ref{btms}.

\subsection{The case $n=2$ for TMS-B}
\label{caseb}

Here we show how the topological mirror symmetry works for $\cM_\B$ when $n=2$. It is instructive to compare it to the arguments in \S\ref{higgs}.\ref{n=2tmsdol}.  

The variant part of the Serre polynomial of $\cM_\B$ is the difference of the full Serre polynomial and the invariant part. This difference in turn can be evaluated using the character tables of $\SL_2(\F_q)$ and $\GL_2(\F_q)$ respectively.
We get that this difference 
  \beq{difference} \nonumber&\sum_{1\neq \kappa\in \hat{\Gamma}}&E_{\kappa}(\cM_\B)=E(\cM_\B)-E_1(\cM_\B) = E(\cM_\B)-E(\calM_\B)/(q-1)^{2g}\\ &=& \nonumber\sum_{\chi \in Irr(\SL_2(\mathbb{F}_q))}\frac{|\SL_2(\mathbb{F}_q)|^{2g-2}}{\chi(1)^{2g-2}}\cdot \frac{\chi(- I)}{\chi(1)}-\sum_{\chi \in Irr(\GL_2(\mathbb{F}_q))}\frac{|\PGL_2(\mathbb{F}_q)|^{2g-2}}{(q-1)\chi(1)^{2g-2}}\cdot \frac{\chi(- I)}{\chi(1)}\\ \nonumber&=& (2^{2g}-1)q^{2g-2}\left( \frac{(q-1)^{2g-2}-(q+1)^{2g-2}}{2}\right)\\&=&\nonumber\sum_{i=1}^{g-1}(2^{2g}-1)
	{\bino{2g-2}{2i-1} q^{2g-3+2i}},\\ \!\! &\, & \ \  \eeq is exactly given  by those $4$ irreducible characters of $\SL_2(\F_q)$ which arise from  irreducible characters of $\GL_2(\F_q)$ which split into two irreducibles over $\SL_2(\F_q)$.
	
		The mapping class group of $C$ acts by automorphisms on  $\Gamma\cong H^1(C,\Z_2)$ 
	so that the induced action on the set $\hat{\Gamma}^*=\hat{\Gamma}\setminus \{1\}$ with $(2^{2g}-1)$ elements is transitive. This way we can argue that $E_{\kappa_1}(\cM_\B)=E_{\kappa_2}(\cM_\B)$ for any two $\kappa_1,\kappa_2\in \hat{\Gamma}^*$, thus we can conclude \beq{one}E_{\kappa}(\cM_\B)=q^{2g-2}\left( \frac{(q-1)^{2g-2}-(q+1)^{2g-2}}{2}\right)\eeq for $\kappa\in \hat{\Gamma}^*$.
	
	  Now $\cM^\gamma_\B$ can be identified with $(\C^\times)^{2g-2}$, which is the Betti version of $\cM^\gamma$ from \S\ref{higgs}.\ref{n=2tmsdol}, which was isomorphic to the cotangent bundle of the identity component of the Prym variety. Now the $\Gamma$-equivariant local system $L_{B,\gamma}$ kills exactly the even cohomology and so we get  \beq{string}E(\cM_\B^\gamma/\Gamma,L_{\hat{B},\gamma})=\frac{(q-1)^{2g-2}-(q+1)^{2g-2}}{2}.\eeq
This proves: \begin{theorem} When $n=2$ Conjecture~\ref{bettirefined} and so Conjecture~\ref{btms} hold. \end{theorem}    As the character tables of $\SL_2({\F_q})$ an	d $\GL_2(\F_q)$ were already known to Schur \cite{schur} and Jordan \cite{jordan} in 1907; in principle the mirror symmetry pattern \eqref{difference} above could have been checked more than 100 years ago! Also the agreement of \eqref{one} and \eqref{string} up to a $q$-power, that is Conjecture~\ref{bettirefined}, when written in terms of character sums  has a somewhat similar form as the fundamental lemma. Indeed in our last section \S\ref{solving}.\ref{TMStoFL} we will argue that Conjecture~\ref{bettirefined} and the fundamental lemma have a common geometrical root. 

\section{Solving our problems}
\label{solving}

Although the arithmetic technique discussed in the previous session  \cite{hausel-etal} is capable of proving the Betti version of the topological mirror symmetry Conjecture~\ref{btms},  it still introduces a different set of conjectures. This raises the question: which one is the "right" one, that is the true consequence of mirror symmetry, the De Rham Conjecture \ref{drtms} and the equivalent Dolbeault Conjecture~\ref{doltms} or the Betti version Conjecture~\ref{btms} of the topological mirror symmetries? 

In the next section we will see that in fact they have a common generalization
which at the same time also solves the earlier "agreement of Hodge numbers" problem discussed in Remark~\ref{handwaving}.

	\subsection{Hard Lefschetz for Weight and Perverse Filtrations}

		We start with the observation \cite[Corollary 3.5.3]{hausel-villegas}  for $\GL_n$ and $\PGL_n$ and \cite{mereb} for $\SL_n$ that	  \beq{palin} E(\hat{\cM}_\B;1/q)=q^{\dim} \ E(\hat{\cM}_\B;q)\eeq  i.e. that the Serre polynomials of our character viarieties are palindromic. (Here $\dim$ is the dimension of the appropriate character variety.) It is interesting to note that ultimately this is due to  Alvis-Curtis duality in the character theory of finite groups of Lie-type. 
	
	We recall  the weight filtration: $$W_0\subset\dots\subset W_{i}\subset\dots \subset W_{2k}=H^k(\M_\B)$$ on the ordinary cohomology $H^*(\M_\B)$. 
		The	 palindromicity \eqref{palin} then lead us to the following Curious Hard Lefschetz Conjecture in \cite[Conjecture 4.2.7]{hausel-villegas}: 
			\beq{curious}\begin{array}{cccc}L^l:& Gr^W_{\dim-2l}H^{i-l}(\M_\B)&\stackrel{\cong}{\rightarrow}& Gr^W_{\dim+2l}H^{i+l}(\M_\B)\\ 
	&x&\mapsto& x\cup \alpha^l		
			\end{array},\eeq  where $\dim=\dim(\M_\B)$ and $\alpha\in W_4H^2(\M_\B)$. For $n=2$ this was proved in \cite[\S5.3]{hausel-villegas}. The conjecture \eqref{curious} is curious because because $\M_\B$ is an affine, thus non-projective, variety and $\alpha$ is a weight $4$ and type $(2,2)$ class, instead of the usual weight $2$ and type $(1,1)$ class of the Hard Lefschetz theorem. However there is a  situation where  a  similar Hard Lefschetz theorem was observed. 
			 
Namely, \cite{decataldo-migliorini} introduce	the 	 perverse filtration: $$P_0\subset \dots \subset P_{ i}\subset \dots P_{k}(X)\cong H^k(X)$$   for $f:X\to Y$ proper $X$ and $Y$ smooth, quasi-projective. Originally they define it using
the BBDG-decomposition theorem of $f_*(\underline{\Q})$ into perverse sheaves. But in \cite{decataldo-migliorini2} they
prove a more elementary equivalent definition in the case when $Y$ is additionally affine. 
Take 	$ Y_0\subset  \dots \subset Y_i \subset \dots Y_d=Y$  s.t. $Y_i$ sufficiently generic with $\dim(Y_i)=i$ then the perverse filtration is given by
					$$P_{k-i-1}H^k(X)={\rm{ker}}(H^k(X)\to H^k(f^{-1}(Y_i))).$$ 
				 Now the Relative Hard Lefschetz Theorem \cite[Theorem 2.3.3]{decataldo-migliorini} holds:
				$$\begin{array}{cccc}L^l:& Gr^P_{{\rm rdim}\!-l}H^*(X)&\stackrel{\cong}{\rightarrow}& Gr^P_{{\rm rdim}+l}H^{*+2l}(X)\\ 
			&x&\mapsto& x\cup \alpha^l		
				\end{array}$$ where ${\rm rdim}$ is the relative dimension of $f$ and $\alpha\in W_2H^2(X)$ is a relative ample class. 
				
					Recall from Theorem~\ref{proper} that the Hitchin map 
			$$\begin{array}{cccc}\chi:&\M_\Dol&\to& \A\\ &(E,\phi)&\mapsto& \rm{charpol}(\phi) \end{array}$$ is proper, thus induces perverse filtration on $H^*(\M_\Dol)$. A nice explanation for the curious Hard Lefschetz 
			conjecture \eqref{curious} would be if the following conjecture held.

			\begin{conjecture}\label{P=W} We have  $P=W$, more precisely, $P_{k}(\M_\Dol)\cong W_{2k}(\M_\B)$ under the  isomorphism $H^*(\M_\Dol)\cong H^*(\M_\B)$ from non-Abelian Hodge theory. \end{conjecture} 
				In \cite{thmhcv} it was proved that
			\begin{theorem}\label{p=w}				$P=W$  when $G=\GL_2,\PGL_2,\SL_2$. \end{theorem}
			\begin{remark} The proof of Theorem~\ref{p=w} was accomplished by a careful study of the topology of the Hitchin map, which paralleled special cases of results of Ng\^{o} in his proof \cite{ngo} of the fundamental lemma.  Additionally we had to use all the previously established results \cite{hauselint,HT1,HT2,hausel-villegas} on the cohomology of these   $n=2$ varieties. Interestingly for the proof for $\SL_2$ we had to use results which were discussed here for the topological mirror symmetry presented in \S\ref{higgs}.\ref{n=2tmsdol} and \S\ref{betti}.\ref{caseb} in this paper. We will now explain why this connection to mirror symmetry is not surprising.
			\end{remark}
	  \subsection{Perverse topological mirror symmetry}
			 We define the perverse Serre polynomial as $$PE(\M_\Dol;u,v,q):=\sum q^k E(Gr^P_k(H_c^*(\M_\Dol));u,v),$$ and the $\hat{B}^d$-twisted stringy perverse Serre polynomial as $$P\Est^{\hat B^d}(\hM^e_\Dol;u,v,q):=\sum_{\gamma\in \Gamma}  PE(\cM_{\Dol,\gamma}^e/\Gamma,L_{\hat{B}^d};u,v,q)(uvq)^{F(\gamma)}.$$
		By Definition \ref{epoly} and Theorem \ref{dr=dol} we have	 \beq{deform}PE(\M_\Dol;u,v,1)=E(\M_\Dol;u,v)=E(\M_{\rm DR};u,v).\eeq
			  Conjecture~\ref{P=W} that $P=W$ then would imply \beq{pespec}PE(\M_\Dol;1,1,q)=E(\M_\B;q)\eeq and Relative Hard Lefschetz \cite[Theorem 2.3.3]{decataldo-migliorini} shows \beq{perversesymmetry}PE(\M_\Dol;u,v,q)=(uvq)^{\dim} PE\left(\M_\Dol;u,v;\frac{1}{quv}\right).\eeq 
			  Note that although the original Hodge numbers of $\M_\Dol$ did not possess any non-trivial symmetry, this refined version with the perverse filtration does. So in fact with this definition we can write down our most general form of the 
			topological mirror symmetry conjectures. 	
				\begin{conjecture} \label{tms}
	$PE\Left(\cM^d_\Dol;x,y,q\Right) = (xyq)^{\dim} P\Est^{\hat B^d}\Left({\hM}^e_\Dol;x,y,\frac{1}{qxy}\Right)$
			\end{conjecture}
\begin{remark}	This most general form of our topological mirror symmetry conjecture solves the two problems we encountered before. 			

First, we see that this version of the topological mirror symmetry conjecture implies
			 Conjecture~\ref{tms} via \eqref{deform} and Conjecture~\ref{btms} if we assume Conjecture~\ref{p=w} and thus \eqref{pespec}. Thus Conjecture~\ref{tms} is a common generalization of Conjectures \ref{drtms}, \ref{doltms} and \ref{btms}. 
			 
Second, relative hard Lefschetz endows the perverse Hodge numbers with the symmetry \eqref{perversesymmetry}, and so one can formulate topological mirror symmetry as in Conjecture~\ref{tms}. On the level of Hodge numbers this conjecture takes the form: 
\beq{hijptms}h^{i,j}_p(\cM_\Dol^d)=h^{i+(\hp-p)/2,j+(\hp-p)/2}_{{\rm st},\hp}(\hM_\Dol^e,\hat{B}^d),\eeq where $\hat{p}=\dim(\cM_\Dol^d) -p$ is the opposite perversity, the Hodge numbers are defined as $$h^{i,j}_p(\cM_\Dol^d):=\dim(H^{i,j}(Gr^W_pH^{i+j}(\cM_\Dol^d)))$$ and similarly for the stringy extension on the right hand side.  This form \eqref{hijptms} is now more reminiscent of the original topological mirror symmetry \eqref{otms}. One can also compare the  functional equation forms \eqref{functional} and Conjecture~\ref{tms}.

The last ingredient which is missing to completely justify our topological mirror symmetry conjectures is to show that indeed these are  cohomological shadows of the S-duality in the work of  Kapustin-Witten \cite{kapustin-witten}. Such an argument will be sketched in the next section.

		\end{remark}

\subsection{Topological mirror symmetry as  cohomological shadow of S-duality}
\label{shadow}
		 Recall that Kapustin and Witten \cite{kapustin-witten} suggest that the Geometrical Langlands program is S-duality  reduced to $2$ dimensions.
		This  simplifies  to  $T$-duality, as first suggested by \cite{vafa-etal}. In turn by the SYZ proposal we get to mirror symmetry between $\cM_\DR$ and $\hM_\DR$. 
		
	 Now Kontsevich's \cite{kontsevich} homological mirror symmetry conjecture implies that \beq{glc}{\mathcal{D}}^b(Coh(\cM_{\DR}))\sim {\mathcal{D}}^b(Fuk(\hM_{\DR})) \eeq the derived category of coherent sheaves on $\cM_\DR$ is equivalent with a certain Fukaya category on $\hM_\DR$. This latter is not straightforward to define but recent work of Nadler-Zaslow  \cite{nadler-zaslow, nadler} relates 
	a certain Fukaya category of $T^*X$ and a category of $D$-modules on $X$, for a compact real analytical manifold $X$. Thus we may imagine that this result
	might extend to give an equivalence of the right hand side of \eqref{glc} 
	with some category of $D$-modules on the stack ${\rm Bun}_{\PGL_n}$ of
	$\PGL_n$ bundles on $C$. 
 	The mathematical content of \cite{kapustin-witten} maybe phrased that the combination of this latter Nadler-Zaslow type equivalence with the homological mirror symmetry in \eqref{glc}  leads to the proposed 
		Geometric Langlands program of \cite{laumon,beilinson-drinfeld}. As explained in \cite{donagi-pantev} in a certain semi-classical
		limit \eqref{glc} should become
		  $${\mathcal{D}}^b(Coh(\cM_{\Dol}))\sim {\mathcal{D}}^b(Coh(\hM_{\Dol})),$$ an equivalence of the derived categories of sheaves on Hitchin systems for Langlands dual groups.  By recent work of Arinkin \cite{arinkin}  it is expected that there is a geometrical  fibrewise Fourier-Mukai transform at least for integral spectral  curves. It means that there should be a Poincar\'e bundle $\calP$ on the fibered product $\cM_\Dol\times_{\A^0}\hM_\Dol$
		 such that the associated fiberwise Fourier-Mukai transform would identify \beq{fourier-mukai} {\mathcal FM}=\hat{\pi}_* ({\calP} \otimes \check{\pi}^*):{\mathcal{D}}^b(Coh(\cM_{\Dol}))\stackrel{\sim}{\to} {\mathcal{D}}^b(Coh(\hM_{\Dol})).\eeq
							One can argue that the cohomological shadow of such a Fourier-Mukai transform for orbifolds should be defined in stringy cohomology. Also if one twists the above Fourier-Mukai transform by adding gerbes, as discussed e.g. in \cite{ben-bassat}, then
							we should see stringy cohomology twisted with gerbes. Also we should expect
							the cohomological shadow of \eqref{fourier-mukai} to be compatible with the perverse filtration. All in all, the cohomological shadow of \eqref{fourier-mukai}	should identify \beq{fmtms}S:H_{p}^{r,s}(\cM^d_{\Dol})\cong H_{st,\hp}^{r+(\hp-p)/2,s+(\hp-p)/2}(\hM^e_\Dol;\hat{B}^d).\eeq In fact, this statement over the regular locus $\A_{\rm reg}$ can be proved. Moreover comparing supports of \eqref{fmtms} over $\A^0$ or also by a Fourier transform argument on $\Gamma$ one can deduce the refined version of \eqref{fmtms}. Namely for $\kappa\in \hat{\Gamma}$ and $\gamma=w(\kappa)\in \Gamma$ we have:
							\beq{fmrefinedtms}S:H_{p}^{r,s}(\cM^d_{\Dol})_\kappa\cong H_{\hp-F(\gamma)}^{r+(\hp-p)/2-F(\gamma),s+(\hp-p)/2-F(\gamma)}(\cM^e_{\Dol,\gamma}/\Gamma;L_{\hat{B}^d}).
							\eeq
							
								 This way we can argue that the cohomological shadow of S-duality reduced to $2$-dimensions and in the semi-classical limit should yield our Topological Mirror Symmetry Conjecture~\ref{tms}:	$$ PE\Left(\cM^d_\Dol;x,y,q\Right) = (xyq)^{\dim} P\Est^{\hat B^d}\Left( {\hM}^e_\Dol;x,y,\frac{1}{qxy} \Right).$$  
		
\begin{remark} 	More delicate structures of the Kapustin--Witten reduced S-duality \cite{kapustin-witten} have been mathematically implemented in recent works of Yun	\cite{yun}, where Ng\^o's techniques from \cite{ngo} were havily used.
\end{remark}				 
								 Finally, in the last section, we explain a connection between our topological mirror symmetry conjectures and Ng\^{o}'s work \cite{ngo} on the fundamental lemma in the Langlands program.

\subsection{From topological mirror symmetry to the fundamental lemma}
\label{TMStoFL}
	
	Ng\^{o}'s celebrated\footnote{A detailed survey of the statement and some of the proof of the fundamental lemma could be found in \cite{nadler2}.}  proof \cite{ngo} of the fundamental lemma is the culmination 
	of a series of  geometrical advances in the understanding of orbital integrals including \cite{kazhdan-lusztig, goresky-etal, laumon2, laumon3, laumon-ngo}. 
	The proof proceeds by studying
	the Hitchin fibration over the so-called elliptic locus. In the case of $\SL_n$, which we will be only discussing here, this means the locus $\A_{\ell}\subset \A^0$ containg characteristics $a\in \A^0$ so that the corresponding spectral curve $C_a$ is integral, i.e. irreducible and reduced. In particular, it contains the locus we studied in this survey $\A_\reg$ where $C_a$ is smooth.
	He considers the degree $0$ Hitchin fibration $\cchi_\ell:\cM^0_{\ell}=\cchi^{-1}(\A_{\ell})\to \A_{\ell}$ over the elliptic locus. 
	An important ingredient in Ng\^{o}'s proof is the BBDG decomposition theorem of  the derived push forward $\cchi_{\ell*}(\underline{\Q})$ of the constant sheaf on $\cM_{\ell}$  into perverse sheaves. He proves
	 the so-called {\em support theorem}, that in certain cases, including the Hitchin fibration, the perverse components of $\cchi_{\ell*}(\underline{\Q})$ are determined by a relatively small open subset of $\A_\ell$. The proof then is achieved by
	checking the geometrical formula \eqref{maingo} below over this small open subset of $\A_\ell$; yielding the statement over the whole $\A_\ell$. 
	
	An important further geometrical insight of the paper \cite{ngo}  is that  $\cchi_{\ell*}(\underline{\Q})$ should be understood with respect to a certain symmetry of the Hitchin fibration, which we already studied in \S\ref{basic}.\ref{symmetries} over $\A_\reg$. Similarly to our definition of the Prym variety $\check{P}_a$ for a smooth spectral curve $C_a$ in \eqref{smoothprym}, we can define the norm map $\Nm_{C_a/C}: \Pic^0(C_a)\to \Pic^0(C)$ and Prym variety $\check{P}_a:=\ker(\Nm_{C_a/C})$ for an integral spectral curve $\pi_a:C_a\to C$ as well. Again similarly to the smooth case one can construct an action of $\check{P}_a$ on $\cM_a:=\cchi^{-1}(a)$ when $a\in 
\A_\ell$. This way we get an action of the group scheme $\check{P}_\ell$ on $\cM_\ell$. 
This symmetry of the Hitchin fibration will induce an action of $\check{P}_\ell$ on $\cchi_{\ell*}(\underline{\Q})$. 
As the group scheme $\check{P}^0_\ell$, the connected component of the identity in $\check{P}_\ell$, acts trivially on cohomology; this action of $\check{P}_\ell$ on $\cchi_{\ell*}(\underline{\Q})$ will factor through an action of the   group scheme of components $\underline{\Gamma}:=\check{P}_\ell/\check{P}^0_\ell$. This turns out to be a finite group scheme with stalk at $a$ agreeing with $\Gamma_a$ the group of components of $\check{P}_a$. 

The finite group scheme $\underline{\Gamma}$  also connects  nicely with
our finite group  $\Gamma=\Pic^0(C)[n]$. Namely one can easily show that for $\gamma\in \Gamma$ the pull back $\pi_a^*(\gamma)\in \check{P}_a$. This way we get a map $$f:\Gamma\to\Gamma_a,$$ which is shown to be surjective in \cite{hausel-pauly}, where the kernel is also explicitly described.  If we now consider a character 
$\kappa\in\hat{\Gamma}_a$, then we get a character $\kappa f\in \hat{\Gamma}$ and a corresponding $\gamma=w(\kappa f)\in \Gamma$. Then the stalk of a sheaf version of the refined S-duality in \eqref{fmrefinedtms} for the $d=0$ case over $\A_\ell$ followed by relative Hard Lefschetz leads to the isomorphism
\beq{maingo} H_{p}^{r,s}(\cM_a)_\kappa\cong H_{p-F(\gamma)}^{r-F(\gamma),s-F(\gamma)}(\cM_{a,\gamma}/\Gamma)\eeq
							This formula, which we derived here from the cohomological shadow of Kapustin-Witten's reduced S-duality, can be identified with the stalk of Ng\^o's main geometric stabilization formula \cite[Theorem 6.4.2]{ngo} in the case of $\SL_n$. As Ng\^{o} argues in \cite{ngo}, when \eqref{maingo} is proved in positive characteristic, and one takes the alternating trace of the Frobenius automorphism on both sides of \eqref{maingo}, then the resulting formula can be seen to imply the fundamental lemma in the Langlands program in the function field case and in turn by Waldspurger's work 
							\cite{waldspurger}  in the number field case.
							
							As explained in \S\ref{solving}.\ref{shadow} the hope is that one can push \eqref{maingo} or more precisely a sheaf version underlying \eqref{fmrefinedtms} from $\A_\ell$ over the whole of the $\SL_n$-Hitchin base $\A^0$. The $\SL_n$-case of the work of   Chaudouard and Laumon \cite{chaudouard-laumon1,chaudouard-laumon2} managed, by extending Ng\^{o}'s techniques, to do this over $\A_\red$ that is managed to prove \eqref{maingo} for reduced, but possibly reducible, spectral curves $C_a$. This way \cite{chaudouard-laumon1,chaudouard-laumon2} lead to	a proof of the so-called weighted fundamental lemma, which again by earlier results of Waldspurger  and others completed the proof of the full endoscopic functoriality principle of Langlands.  For us however it remains to extend \eqref{fmrefinedtms} over the whole of $\A^0$, including non-reduced and reducible spectral curves, which will yield our
						 topological mirror symmetry Conjecture~\ref{tms}. 
							
							Details of the arguments in the last two sections \S\ref{solving}.\ref{shadow} and \S\ref{solving}.\ref{TMStoFL} will appear elsewhere.

			\section{Conclusion}\label{conclusion}
	 
	 In this paper we surveyed some techniques to obtain cohomological information on the topology of the total space of the Hitchin system. We painted a picture where ideas from physics and number theory were  combined into  a dynamic mix. Although these techniques are fairly powerful, still they have not yet lead to complete understanding. In particular, the most general conjectures are still open.
	 
	 More recently there have appeared work by physicists Diaconescu et al. \cite{diaconescu-etal} about a new string theory framework for several conjectures relating to the topology of the total space of the $\GL_n$-Hitchin system. Besides the links
	 to the conjectured formulae in \cite{hausel-villegas, hausel-mln} a picture is emerging which relates the main conjecture of \cite{hausel-villegas, hausel-mln} with a certain version of the Gopukamar-Vafa conjecture, which ultimately can be phrased as strong support of our pivotal $P=W$ Conjecture~\ref{P=W}. This way \cite{diaconescu-etal} uncovers close connections  of our conjectures  with Gromov-Witten, Donaldson-Thomas and Pandharipande-Thomas invariants of certain local Calabi-Yau $3$-folds. There have been considerable progress on the latter invariants lately in the mathematics literature thus we can well hope that with this new point of view we will be able to progress our understanding of the problems surveyed in this paper.
	 
	 We finish by mentioning another promising new work  \cite{heinloth-etal} where they manage to extend the original Morse theory method
	 of Hitchin for $n=2$ and Gothen for $n=3$ to higher $n$ using a motivic view point - originating in the number theoretic approach of Harder-Narasimhan \cite{harder-narasimhan}. In particular, their calculations have been done for $n=4$ which are in agreement with the conjectures in \cite{hausel-villegas,hausel-mln} and this paper. 
	 
	 When we add these two very recent approaches, again one originating in physics and one in number theory, to the mix of ideas surveyed in this paper, we can be sure that new exciting results and ideas will be found on questions relating to the global topology of the Hitchin system in the foreseeable future.

\end{document}